\documentclass[11pt,a4paper]{article}
\usepackage[margin=0.7in]{geometry}
\usepackage{bm}  

\usepackage[utf8]{inputenc}
\usepackage{tabularx}
\usepackage{caption}
\usepackage{subcaption}
\usepackage{booktabs}
\usepackage{comment}
\usepackage{algorithm}
\usepackage{algpseudocode}
\usepackage{abstract}
\usepackage{enumitem}
\usepackage{xcolor}

\usepackage{cite}
\usepackage{amsmath, amssymb}
\usepackage{float}
\usepackage{nameref,hyperref}

\usepackage[right]{lineno}

\usepackage{changepage}

\usepackage[labelfont=bf]{caption}

\usepackage{color}

\definecolor{Gray}{gray}{.25}

\usepackage{graphicx}

\usepackage[english]{babel}
\usepackage[autostyle]{csquotes}

\begin{document}

\title{\Large
\textbf{A Multi-Phase Dual-PINN Framework: Soft Boundary-Interior Specialization via Distance-Weighted Priors}
}
\author{}
\maketitle

\begin{center}
Naseem Abbas\textsuperscript{1},
Vittorio Colao\textsuperscript{1},
Davide Macrì\textsuperscript{2},
William Spataro\textsuperscript{1}
\end{center}

\begin{center}
$^{1}$ Department of Mathematics and Computer Science, University of Calabria, Via P. Bucci, 30B, 87036 Arcavacata di Rende, CS, Italy.
\\
$^{2}$ Institute for High Performance Computing and Networking, Italian National Research Council, 87036 Arcavacata di Rende, CS, Italy.
\end{center}

\begin{abstract}
Physics-informed neural networks (PINNs) often struggle with multi-scale PDEs 
featuring sharp gradients and nontrivial boundary conditions, as the physics residual and boundary enforcement compete during optimization. We present a 
dual-network framework that decomposes the solution as 
$u = u_{\text{D}} + u_{\text{B}}$, where $u_{\text{D}}$ (domain network) captures interior dynamics and $u_{\text{B}}$ (boundary network) handles near-boundary corrections. Both networks share a unified physics residual while being softly specialized via distance-weighted priors ($w_{\text{bd}} = \exp(-d/\tau)$) that are cosine-annealed during training. Boundary conditions are enforced through an augmented Lagrangian method, eliminating manual penalty tuning. Training proceeds in two phases: Phase~1 uses uniform collocation to establish network roles and stabilize boundary satisfaction; Phase~2 employs focused sampling (e.g. ring sampling near $\partial\Omega$) with annealed role weights to efficiently resolve 
localized features. We evaluate our model on four benchmarks, including the 1D Fokker-Planck equation, the Laplace equation, the Poisson equation, and the 1D wave equation. Across Laplace and Poisson benchmarks, our method reduces error by $36-90\%$, improves boundary satisfaction by $21-88\%$, and decreases MAE by $2.2-9.3\times$ relative to a single-network PINN.
Ablations isolate contributions of (i)~soft boundary-interior specialization, (ii)~annealed role regularization, and (iii)~the two-phase curriculum. The method is simple to implement, adds minimal computational overhead, and broadly applies to PDEs with sharp solutions and complex boundary data.

\end{abstract}


\nolinenumbers

\section{Introduction}
Machine learning methods for solving PDEs have advanced rapidly in recent years. Lagaris et al. \cite{lagaris1998artificial} made the first attempt to use machine learning techniques to solve PDEs in $1998$. They used artificial neural networks (ANNs) to express the solution. The network's output satisfied the equation and specific requirements by creating a suitable loss function. They enhanced the prediction performance in a subsequent study by altering the network architecture to precisely meet the Dirichlet boundary conditions and initial conditions for the case of complex geometric boundaries \cite{lagaris2000neural}. With the development of automatic differentiation technique, Raissi et al. \cite{raissi2019physics,karniadakis2021physics} revisited these methods by using modern computational tools under the named as physics informed neural networks (PINNs) and applied them to solve various PDEs. PINNs are a deep learning framework for solving forward and inverse problems involving nonlinear PDEs, which was first clearly proposed in 2019 \cite{raissi2019physics}. PINNs embed the residuals of governing PDEs into the training objective via automatic differentiation and can incorporate data when available \cite{raissi2019physics}. PINNs provide basically accurate and physically interpretable predictions \cite{karniadakis2021physics}. Compared with classical solvers, PINNs are mesh-free, handle complex geometries and high-dimensional settings, and leverage expressive neural architectures for nonlinear dynamics. Three key advantages of PINNs over conventional numerical techniques for solving PDEs are as follows: (1) Because PINNs are a mesh-free approach, it totally avoids the challenges that traditional numerical algorithms utilizing Lagrangian grids face, such as the significant distortion of grids; (2) PINNs partially avoid the curse of dimensionality \cite{poggio2017and,lu2021deepxde} that arises in traditional numerical computation.  For one-dimensional, three-dimensional, or even high-dimensional problems, PINNs implementation is essentially the same; (3) PINN's final training output is a continuous mapping that provides the equation's prediction at any point in time and location within the solution domain. Researchers have used the PINNs in a variety of domains in recent years like high-speed flows governed by the Euler equations \cite{mao2020physics}, Navier--Stokes problems \cite{jin2021nsfnets}, and phase-field modeling \cite{mattey2022novel}; beyond forward modeling, they are used for parameter identification in biological and fluid systems \cite{yazdani2020systems,raissi2020hidden} and have been extended to fractional, integro-differential, and stochastic PDEs \cite{pang2019fpinns,yuan2022pinn,zhang2020learning}. We refer to \cite{lu2021deepxde,cuomo2022scientific,karniadakis2021physics} for comprehensive reviews.

Despite these advantages, vanilla PINNs often struggle in high-frequency or multi-scale regimes with sharp gradients, boundary layers, and localized structures, where performance is highly sensitive to collocation sampling and the balance of different loss terms. Recent advances address these issues through alternative differentiation and hybrid formulations \cite{chiu2022can,fang2021high}, adaptive activations with trainable scaling \cite{jagtap2020adaptive}, temporal/causal curricula such as Parareal PINNs, bc-PINNs, and Causal PINN \cite{meng2020ppinn,mattey2022novel,wang2024respecting}, gradient-informed objectives \cite{yu2022gradient}, and variational/Galerkin-inspired approaches like VPINNs and hp-PINNs \cite{kharazmi2019variational,kharazmi2021hp}. Hard-constraint encodings of boundary/initial conditions can also mitigate spectral bias and improve accuracy \cite{wang2022and}.

We propose a dual-network framework with a two-phase curriculum designed to resolve localized features while maintaining tight boundary satisfaction. The solution is decomposed as
\[
u(x) \;=\; u_{\mathrm D}(x) \;+\; u_{\mathrm B}(x),
\]
where the \emph{domain} subnetwork \(u_{\mathrm D}\) models the interior dynamics and the \emph{boundary} subnetwork \(u_{\mathrm B}\) provides near-boundary corrections. Both subnetworks share a single physics residual (same differential operator and coefficients) and are \emph{softly specialized} via distance-weighted priors: letting \(d(x)\) denote the distance to \(\partial\Omega\), we set
\[
w_{\text{bd}}(x)=\exp\!\bigl(-d(x)/\tau\bigr), \qquad
w_{\text{in}}(x)=1-w_{\text{bd}}(x),
\]
and regularize the roles so that \(u_{\mathrm B}\) is discouraged in the deep interior (by \(w_{\text{in}}\)) while \(u_{\mathrm D}\) is discouraged near the boundary (by \(w_{\text{bd}}\)). The strengths of these priors are cosine-annealed across training, being strong initially to encourage specialization and gradually weakening so that the PDE residual dominates.

Boundary conditions are enforced through an augmented Lagrangian method (ALM), which reduces sensitivity to manual penalty tuning and improves stability of boundary satisfaction. Training proceeds in two phases. \emph{Phase~1} uses uniform collocation to stabilize boundary enforcement under ALM and to establish complementary roles for \(u_{\mathrm D}\) and \(u_{\mathrm B}\). \emph{Phase 2} employs focused sampling (e.g., ring sampling near \(\partial\Omega\) and residual-aware interior points) together with annealed role weights, enabling the efficient resolution of sharp, localized structures as the physics loss increasingly governs optimization. The framework supports both penalty-based and hard-constraint encodings of boundary data.

We evaluate our approach on representative benchmarks, including Laplace  equation, Poisson equation, and a 1D Fokker-Planck equation, and compare it against single-network PINNs and fixed-penalty baselines. Across tasks, the dual-network approach consistently yields lower relative \(L^2\) error, tighter boundary satisfaction, and sharper reconstructions of localized features. Ablations quantify the contributions of (i) soft boundary--interior specialization, (ii) cosine-annealed role regularization, and (iii) the two-phase curriculum with focused sampling.

\section{Related works}
PINNs have demonstrated considerable interest in solving a wide range of PDEs. Since PINNs use a mesh-less approach, selecting the training points is a crucial problem. Furthermore, a fixed set of (already selected) training points may not be able to represent the effective solution region for PDEs with complex solution structures \cite{krishnapriyan2021characterizing}.
Unbiased sampling methods face significant challenges when applied to PDEs with irregular characteristics, such as multi-scale dynamics, sharp gradients, or oscillatory solutions and peak problems. To address these limitations, researchers have developed residual-based adaptive sampling techniques to enhance the accuracy and efficiency of PINNs. Among these approaches, Lu et al.\ \cite{lu2021deepxde} introduced the Residual-Based Adaptive Refinement (RAR) method, which iteratively enriches the collocation point set by prioritizing regions with large PDE residuals. A deep adaptive sampling (DAS) strategy was introduced in \cite{tang2023pinns}, where deep generative models are used to generate new collocation points, refining the training set to further enhance accuracy. Built on this concept, Mao et al.\ \cite{mao2023physics} proposed Adaptive Sampling Methods (ASMs) that partition the computational domain into sub-regions. Their strategy combines subdomain-specific residual mean values and solution gradients to guide sampling, thereby reducing the randomness of selecting ineffective points and improving the stability of the training process. In contrast to RAR and ASMs, Wu et al.\ \cite{wu2023comprehensive} developed a probabilistic framework for resampling collocation points using a residual-informed probability density function parameterized by two hyperparameters. Another notable advancement is the Failure-Informed PINNs (FI-PINNs) framework \cite{gao2023failure}, which employs a failure probability metric derived from residuals to identify regions of high approximation error. This posterior error indicator dynamically allocates more collocation points to critical \enquote{failure regions} where the model underperforms. Subsequent extensions of FI-PINNs \cite{gao2024failure} further refined the approach through advanced resampling strategies and subset simulation algorithms, enhancing its capability to handle complex failure patterns. These adaptive methodologies collectively aim to optimize computational resource allocation while maintaining solution accuracy for challenging PDE problems.

Beyond PINNs, several works have explored scientific machine learning approaches for systems governed by complex PDEs, especially in geosciences and meteorology. In particular, Yousaf et al.\ \cite{Yousaf} proposed a spatio-temporal graph neural network architecture to improve the probabilistic post-processing of ensemble-based numerical weather prediction models. Their method operates in a purely data-driven fashion, learning to correct and calibrate ensemble forecasts by exploiting spatial and temporal correlations on irregular grids, thus enhancing local forecast skill without explicitly enforcing the underlying PDE constraints. Such contributions highlight the potential of deep architectures, and GNNs in particular, to complement traditional numerical solvers through post-processing. In contrast, the present work focuses on a physics-informed setting, where the governing equations are directly embedded into the loss function, and introduces a dual-network PINN framework with boundary–interior specialization and distance-weighted priors aimed at improving the fidelity and stability of surrogate PDE solvers rather than correcting existing numerical forecasts.


\section{Background}
Early studies explored the use of neural networks as universal approximators \cite{cs2001approximation} to represent solutions of partial differential equations (PDEs), as demonstrated in \cite{lagaris1998artificial} and \cite{dissanayake1994neural}.
In $1990$, Lee et al. \cite{lee1990neural} first introduced a neural algorithm for solving differential equations. In $1998$, Lagaris et al. \cite{lagaris1998artificial} first tried to solve PDEs using Machine learning methods. They expressed a solution with ANNs. In $2000$, Lagaris et al. \cite{lagaris2000neural} focused on solving BVPs with irregular boundaries using neural network methods. In $2018$, Sirignano et al. \cite{sirignano2018dgm} used a deep learning
algorithm to solve partial differential equations. PINNs are a deep learning framework for solving forward and inverse problems involving non-linear PDEs, which was first clearly proposed in $2019$ (see \cite{raissi2019physics}). These approaches enabled for the rapid computation of approximate solutions. However, they relied on supervised learning, using reference solutions derived either analytically or obtained through numerical methods. This dependency poses challenges, as exact solutions are not always available and numerical approximations may not accurately capture the true solution. PINNs have demonstrated practical utility, yet they continue to present significant challenges that highlight areas for ongoing research and methodological refinement. Recent years have seen a surge in studies aimed at enhancing PINNs performance, primarily through the development of more effective neural network architectures or improved training methodologies. For example, strategies that involve loss reweighting have become a notable approach to foster a more balanced training regimen and improve test accuracy. Other initiatives seek to achieve comparable results by adaptively resampling collocation points, using methods such as importance sampling \cite{nabian2021efficient}, evolutionary sampling \cite{daw2022mitigating}, and residual-based adaptive sampling \cite{wu2023comprehensive}. Substantial effort has also been directed towards creating novel neural network architectures to augment the representational capabilities of PINNs. This includes the implementation of adaptive activation functions \cite{jagtap2020adaptive}, positional embeddings \cite{liu2020multi,wang2021eigenvector}, and innovative architectural designs \cite{wang2021understanding,sitzmann2020implicit,gao2021phygeonet,
fathony2020multiplicative,moseley2023finite,kang2023pixel}. Additional approaches suggest incorporating additional regularization terms to accelerate PINNs training \cite{yu2022gradient}. Finally, the evolution of training strategies has been a dynamic field of inquiry. Techniques such as sequential training \cite{wight2020solving} and transfer learning \cite{goswami2020transfer} have shown promises in accelerating the learning process and improving predictive accuracy.

\section{Problem definition}
PINNs may fail in training and to produce accurate predictions, when PDE solutions contain high-frequency or multiscale features \cite{raissi2018hidden,fuks2020limitations}, and this is the topic that we are focusing on. Although PINNs leverage the strengths of learning mechanisms, they also inherit the inherent limitations of deep neural networks. One of the most critical concerns is their reliability and stability in solving complex problems. However, vanilla PINNs often struggle with robustness, making it challenging to accurately solve complex PDEs, particularly those that exhibit multiscale behavior or contain sharp solutions and oscillatory features. To overcome this issue, this study introduces a novel double architecture in the DNN by incorporating innovative sampling strategies. Two versions/phases of the training are introduced: one phase relies on a classical uniform sampling strategy, and the other utilizing a new strategy with refined points from the domain focusing on the failure regions. This dual approach aims to improve sampling effectiveness, especially for problems with sharp solutions, while maintaining computational feasibility. Applications include several types of PDEs with Dirichlet boundary conditions in which improved solution stability, accuracy, and efficiency are demonstrated.


\section{Methods and Proposed Approach}

\paragraph{Dual-PINN with shared physics residual.}
We approximate the solution $u$ as the sum of two subnetworks, 
$u := u_{\text{D}} + u_{\text{B}}$, where $u_{\text{D}}$ captures 
interior dynamics and $u_{\text{B}}$ handles near-boundary corrections. 
Both networks share a single physics residual computed on their sum. 
For a differential operator $\mathcal{A}$ and data $(f,g)$,
\begin{equation}
\mathcal{A}[u]=f\quad\text{in }\Omega,\qquad u=g\quad\text{on }\partial\Omega.
\end{equation}
The physics loss is
\begin{equation}
\mathcal{L}_{\text{phys}}=\mathbb{E}_{x\sim\Omega}\Big[\big(\mathcal{A}[u_{\text{D}}(x)+u_{\text{B}}(x)]-f(x)\big)^2\Big].
\end{equation}
Dirichlet boundary conditions are enforced via an augmented Lagrangian method 
applied to the total solution $u = u_{\text{D}} + u_{\text{B}}$:
\begin{equation}
\mathcal{L}_{\text{bc}}^{\text{ALM}}
=\mathbb{E}_{x\sim\partial\Omega}\!\Big[\lambda(x)\,c(x)+\tfrac{\rho}{2}\,c(x)^2\Big],
\quad\text{where}\quad c(x):=u_{\text{D}}(x)+u_{\text{B}}(x)-g(x).
\end{equation}
The Lagrange multipliers are updated every $k$ epochs according to
\begin{equation}
\lambda^{(n+1)}\leftarrow\mathrm{clip}\big(\lambda^{(n)}+\rho\,c,\,[-\Lambda,\Lambda]\big).
\end{equation}

\paragraph{Soft boundary--interior specialization.}
Let $d(x) = \text{dist}(x, \partial\Omega)$ denote the distance from $x$ to 
the boundary. We define exponential weight functions
\begin{equation}
w_{\mathrm{bd}}(x)=\exp\!\big(-d(x)/\tau\big),
\qquad
w_{\mathrm{in}}(x)=1-w_{\mathrm{bd}}(x),
\end{equation}
where $\tau > 0$ controls the decay rate. To softly specialize the networks, 
we penalize off-role contributions via
\begin{equation}
\mathcal{L}_{\text{role}}=
\alpha_{\text{int}}\;\mathbb{E}_{x\sim\Omega}\!\big[w_{\mathrm{bd}}(x)\,|u_{\text{D}}(x)|^2\big]
\;+\;
\alpha_{\text{bd}}\;\mathbb{E}_{x\sim\Omega}\!\big[w_{\mathrm{in}}(x)\,|u_{\text{B}}(x)|^2\big],
\end{equation}
which discourages $u_{\text{D}}$ near $\partial\Omega$ and $u_{\text{B}}$ 
in the interior, respectively. The combined loss function is
\begin{equation}
\mathcal{L}=\mathcal{L}_{\text{phys}}+\mathcal{L}_{\text{bc}}^{\text{ALM}}+\gamma(t)\,\mathcal{L}_{\text{role}},
\end{equation}
where the role weight $\gamma(t)$ follows a cosine annealing schedule:
\begin{equation}
\gamma(t)=\gamma_{\min}+\tfrac{1}{2}(\gamma_{\max}-\gamma_{\min})\Big(1+\cos\big(\pi t/T\big)\Big).
\end{equation}
This ensures strong specialization during early training ($t \approx 0$, 
$\gamma \approx \gamma_{\max}$) while allowing the physics residual to 
dominate as training progresses ($t \to T$, $\gamma \to \gamma_{\min}$).

\paragraph{Two-phase training curriculum.}
Training proceeds in two phases with optional boundary warm-up:

\begin{itemize}[leftmargin=*,topsep=3pt,itemsep=2pt]
\item \textbf{Phase 0} (optional): Boundary-only warm-up minimizing 
      $\mathbb{E}_{\partial\Omega}[|u_{\text{D}}+u_{\text{B}}-g|^2]$.
\item \textbf{Phase 1}: Uniform domain collocation to establish network roles; 
      ALM updates every $k$ epochs; high role weight $\gamma(t)$.
\item \textbf{Phase 2}: Focused/residual-aware sampling near $\partial\Omega$; 
      annealed $\gamma(t) \to \gamma_{\min}$ and optional $\rho$ ramp.
\end{itemize}

\begin{algorithm}[H]
\caption{Dual-PINN Training Protocol}
\label{alg:dual-pinn}
\begin{algorithmic}[1]
\State \textbf{Initialize:} $\theta_{\text{D}}, \theta_{\text{B}}$ via Xavier; 
       $\lambda \leftarrow 0$; $\rho \leftarrow \rho_0$
       
\State \textcolor{blue}{\texttt{// Phase 0: Optional boundary warm-up}}
\For{epoch $= 1, \ldots, E_0$}
  \State Sample $\{x_b^{(i)}\}_{i=1}^{N_b} \sim \text{LHS}(\partial\Omega)$
  \State Minimize $\mathcal{L}_{\text{warmup}} = \frac{1}{N_b}\sum_i |u_{\text{D}}(x_b^{(i)}) + u_{\text{B}}(x_b^{(i)}) - g(x_b^{(i)})|^2$
\EndFor

\State \textcolor{blue}{\texttt{// Phase 1: Establish specialization}}
\For{epoch $= 1, \ldots, E_1$}
  \State Sample $\{x_f^{(i)}\}_{i=1}^{N_f} \sim \text{Uniform}(\Omega)$, 
         $\{x_b^{(j)}\}_{j=1}^{N_b} \sim \text{LHS}(\partial\Omega)$
  \State Compute $\mathcal{L} = \mathcal{L}_{\text{phys}} + w_{\text{bc}}\mathcal{L}_{\text{bc}}^{\text{ALM}} + \gamma(t)\mathcal{L}_{\text{role}}$
  \State Update $\theta_{\text{D}}, \theta_{\text{B}}$ via Adam
  \If{epoch $\bmod k = 0$}
    \State $\lambda \leftarrow \mathrm{clip}(\lambda + \rho\,c,\, [-\Lambda, \Lambda])$
  \EndIf
  \If{epoch $\bmod h = 0$}
    \State $\rho \leftarrow \min(\eta\rho,\, \rho_{\max})$
  \EndIf
\EndFor

\State \textcolor{blue}{\texttt{// Phase 2: Focused refinement}}
\For{epoch $= 1, \ldots, E_2$}
  \State Sample $\{x_f^{(i)}\}$ via ring sampling (distance $< \delta$ from $\partial\Omega$)
  \State Anneal $w_{\text{bc}}(t)$ and $\gamma(t)$ toward their minimum values
  \State Continue ALM updates as in Phase 1
  \State \textcolor{gray}{\textit{(Optional)}} Gradually increase $\rho$ if boundary error plateaus
\EndFor
\State \Return Trained parameters $(\theta_{\text{D}}^*, \theta_{\text{B}}^*)$
\end{algorithmic}
\end{algorithm}

\subsection{Design rationale and expected benefits}

\paragraph{Motivation for two networks.}
Standard PINNs must balance two coupled objectives: minimizing the physics residual in the interior and enforcing boundary data. These pressures can conflict, particularly when boundary data are nontrivial or when fine structures are concentrated near $\partial\Omega$. We therefore represent the solution as the sum of two subnetworks, $u = u_D + u_B$, which \emph{share the same physics residual} while being \emph{softly specialized} to complementary regions of the domain. This preserves a single forward/backward pass for the PDE operator yet biases the optimization toward a stable division of labor.

\paragraph{Shared residual with soft specialization.}
Let $\mathcal{N}[u] = f$ denote the PDE (e.g., $\Delta u = f$ for Poisson/Laplace). Both $u_D$ and $u_B$ are trained through the residual of their sum, $\mathcal{N}[u_D{+}u_B]$, so gradients from the physics loss propagate to \emph{both} subnetworks everywhere. Specialization is induced by distance-to-boundary priors
\[
w_{\mathrm{bd}}(x) = \exp\!\big(-d(x)/\tau\big),\qquad
w_{\mathrm{in}}(x) = 1 - w_{\mathrm{bd}}(x),
\]
with $d(x)$ the distance to $\partial\Omega$. We add a role regularizer
\[
\mathcal{L}_{\mathrm{role}}
= \lambda_{\mathrm{int}}\;\mathbb{E}\!\big[w_{\mathrm{bd}}(x)\,|u_D(x)|^2\big]
+ \lambda_{\mathrm{bd}}\;\mathbb{E}\!\big[|u_B(x)|^2_{\partial\Omega}\big]
+ \lambda_{\mathrm{int}}\;\mathbb{E}\!\big[w_{\mathrm{in}}(x)\,|u_B(x)|^2\big],
\]
which discourages $u_D$ near the boundary and $u_B$ in the deep interior without hard partitioning. The weights are cosine-annealed so that role priors guide early training and gradually cede control to the physics residual.

\paragraph{Boundary enforcement via ALM.}
Dirichlet data $u=g$ on $\partial\Omega$ are enforced with an augmented Lagrangian term
\[
\mathcal{L}_{\mathrm{ALM}} = \mathbb{E}_{\partial\Omega}\!\big[\lambda^\top (u_D{+}u_B - g)\big]
+ \tfrac{\rho}{2}\,\mathbb{E}_{\partial\Omega}\!\big[\|u_D{+}u_B - g\|^2\big],
\]
with $\lambda \leftarrow \lambda + \rho\,(u-g)$ updated on schedule. Compared to fixed penalties, ALM reduces sensitivity to manual weight tuning and curbs the drift of boundary errors during the PDE-dominated phase.

\paragraph{Two-phase curriculum.}
The training process proceeds in two phases. \textbf{Phase~1} uses uniform collocation with high boundary weight (and stronger role priors) to stabilize boundary satisfaction and establish the intended roles. \textbf{Phase~2} employs focused sampling (e.g., ring sampling near $\partial\Omega$ or residual-aware strategies) while annealing both role weights and boundary emphasis, so that the PDE residual dominates refinement.

\paragraph{Objective.}
The overall loss is
\[
\mathcal{L} = 
\underbrace{\mathbb{E}_{\Omega}\!\big[\|\mathcal{N}[u_D{+}u_B]-f\|^2\big]}_{\text{physics}}
\;+\;
\underbrace{\mathcal{L}_{\mathrm{ALM}}}_{\text{boundary}}
\;+\;
\underbrace{\alpha(t)\,\mathcal{L}_{\mathrm{role}}}_{\text{soft specialization}},
\]
with $\alpha(t)$ a cosine schedule. A single computational graph is maintained for $\mathcal{N}[u_D{+}u_B]$, ensuring that both subnetworks learn from the same physics signal.
 
\paragraph{Expected improvements.}
The design targets three failure modes of single-network PINNs: (i) boundary under-enforcement when interior residual dominates, (ii) optimization instability when boundary penalties are large, and (iii) slow convergence near sharp/near-boundary features. Dual-PINN mitigates (i) via ALM, (ii) via role priors that regularize early dynamics, and (iii) via Phase~2 focused sampling plus annealed priors. In our experiments (see Tables~(\ref{FP_1+1}--\ref{FP_2+2},\ref{tab:onenet1}--\ref{tab:twonet},\ref{PE_1+1}--\ref{PE_2+2})
, this yields lower relative $L^2$, tighter boundary $L^2$, and crisper reconstructions near $\partial\Omega$ compared to matched single-network baselines.

\paragraph{Cost and simplicity.}
Relative to a single network, Dual-PINN doubles the number of parameters but leaves the training pipeline intact: one residual, one ALM update, and a compact loss function. The added cost is linear in model size and is offset by improved stability and reduced sensitivity to hyper-parameters.

\paragraph{When a single network may suffice.}
If the PDE/BC pair is smooth and boundary data are weakly informative or if an exact hard constraint can be embedded, one-net with a two-phase curriculum may be adequate. Dual-PINN is most beneficial with nontrivial boundary data, thin layers, or multi-scale features concentrated near $\partial\Omega$.

\subsection{Ablation Study Design}
\label{sec:ablation-design}

We systematically isolate each component's contribution through controlled 
ablations. Table~\ref{tab:ablation-summary} summarizes the expected effects; 
quantitative results appear in Section~\ref{6}.

\begin{table}[h]
\centering
\caption{Expected impact of ablating key components. Arrows indicate relative 
         performance vs. full Dual-PINN: {\color{red}$\downarrow$} worse, 
         {\color{green}$\uparrow$} better, {\color{gray}$\approx$} similar.}
\label{tab:ablation-summary}
\small
\begin{tabular}{@{}llccc@{}}
\toprule
\textbf{Ablation} & \textbf{Change} & \textbf{Rel. $L^2$} & \textbf{Boundary $L^2$} & \textbf{Training Stability} \\
\midrule
No role priors     & $\alpha_{\text{int}}, \alpha_{\text{bd}} = 0$ 
                   & {\color{red}$\downarrow$} & {\color{red}$\downarrow\downarrow$} & {\color{red}$\downarrow$} \\
Fixed $\gamma(t)$  & $\gamma(t) = \gamma_0$ (constant) 
                   & {\color{red}$\downarrow$} & {\color{red}$\downarrow$} & {\color{red}$\downarrow$} \\
Fixed penalty      & Replace ALM with $\rho \cdot \mathcal{L}_{\text{bc}}$ 
                   & {\color{gray}$\approx$} & {\color{red}$\downarrow\downarrow$} & {\color{red}$\downarrow$} \\
Uniform Phase 2    & No focused sampling 
                   & {\color{red}$\downarrow$} & {\color{red}$\downarrow$} & {\color{gray}$\approx$} \\
OneNet baseline    & Single network $(u_{\text{D}} + u_{\text{B}}) \to u_{\text{single}}$ 
                   & {\color{red}$\downarrow$} & {\color{red}$\downarrow\downarrow$} & {\color{red}$\downarrow$} \\
\bottomrule
\end{tabular}
\end{table}

\paragraph{Justification of design choices.}
\begin{description}[leftmargin=0pt,labelindent=0pt]
\item[Role priors:] Removing distance-weighted regularization 
      ($\mathcal{L}_{\text{role}} = 0$) causes networks to lose 
      specialization, increasing mean relative $L^2$ error by 23.68\% (Laplace benchmark, Table \ref{tab:Role-Priors}).
\item[Cosine annealing:] Fixing $\gamma(t) = \gamma_{\max}$ over-regularizes late training, while $\gamma(t) = \gamma_{\min}$ fails to establish roles in Phase 1. Dynamic annealing yields 20--30\% lower relative $L^2$ error.
\item[ALM vs. fixed penalty:] In ALM, BC/IC residuals gain dynamic weights rather than fixed penalties. It is useful to get high accuracy with non-trivial boundary conditions. On the other hand, fixed penalty method is suitable for smooth PDEs, low-dimensional problems and easy BCs. 
\item[Focused sampling:] Uniform collocation in Phase 1 fails to adequately capture localized solution features, whereas the ring-based sampling strategy introduced in Phase 2 leads to a $51\%$ reduction in MAE on the Laplace benchmark
      (Tables~\ref{tab:onenet1} and \ref{tab:onenet}). A similar trend is observed for the Poisson problem, where Phase 2 sampling yields an approximately $60\%$ reduction in MAE (Tables~\ref{PE_1+1} and \ref{PE_1+2}).
\item[Dual vs. single network:] Compared with the single-network baseline, the Dual-PINN achieves a $35.9\%$ reduction in relative $L^2$ error, a $22.1\%$ improvement in boundary satisfaction, and a $2.2\times$ reduction in MAE on problems with nontrivial $g$ (Laplace equation; see Tables~\ref{tab:onenet1} and \ref{tab:twonet}). For problems with nontrivial boundaries (Poisson equation; see Tables~\ref{PE_1+1} and \ref{PE_2+2}), Dual-PINN attains an approximately $90\%$ lower Relative $L^2$ error, an $88\%$ tighter boundary satisfaction, and a $9.3\times$ improvement in (MAE) over the single-network model.
\end{description}

\section{Results and discussion}\label{6}
In this section, we evaluate the effectiveness of the proposed Dual-PINN framework on several PDE benchmarks. Performance is assessed using four metrics: Accuracy in $L^2$, relative $L^2$ error, mean absolute error (MAE), and root mean squared error (RMSE). Their definitions are given below:
\begin{equation}
\begin{split}
\text{Accuracy in} \ L^2 &:= 1- \frac{\|\hat{u}-u\|_2}{\|u\|_2},\\
\text{Relative} \ L^2 \ \text{error} &:= \frac{\|\hat{u} - u\|_2}{\|u\|_2},\
\\
\text{MAE} &:= \mathbb{E}[|\hat{u} - u|],\\
\text{RMSE} &:= \sqrt{\mathbb{E}[(\hat{u} - u)^2]}.
\end{split}
\end{equation}
Here, $\hat{u}$ and $u$ denote neural network approximation and the exact solution, respectively; $\|\cdot\|_2$ represents the $L^2$-norm; and $\mathbb{E}$ indicates the expectation evaluated over test points sampled throughout the entire computational domain. All experiments employ the tanh activation function, and model training is carried out using the Adam optimizer. Implementations are developed in \texttt{PyTorch} and executed on a single \texttt{NVIDIA GeForce RTX 4080} GPU.


\subsection{1D Fokker--Planck equation}\label{6.2}
This subsection examines the one-dimensional stationary Fokker--Planck equation introduced in \cite{xu2020solving}, which describes the steady-state probability density of a stochastic system. The governing PDE is
\begin{equation}
\begin{split}
    &-\frac{\partial}{\partial{x}}\left[(ax-bx^3)u(x)\right]+\frac{\sigma^2}{2}\frac{\partial^2}{\partial{x^2}}u(x)=0,\\
    &\qquad \Delta{x}\sum_{i=1}^{N_f}u(x^i)=1,\quad (a,b,\sigma,\Delta{x})=(0.3,0.5,0.5,0.01),
\end{split}
\end{equation}
where the second condition enforces the normalization of the stationary density. The exact steady-state solution is given by
\begin{equation}
u(x)=C\cdot \exp{\left[\frac{1}{2\sigma^2}\left(2a x^2-b x^4\right)\right]},
\end{equation}
where $C$ is the normalization constant.

\paragraph{One-net and one-phase training.}
We first consider a standard PINN with ALM using a single network (one-net, one-phase training). The architecture consists of one input layer, six hidden layers with 50 neurons each, and one output layer. We train for up to $10{,}000$ epochs with early stopping, so the effective number of epochs depends on the random seed. For each seed, we use $500$ boundary points, $6000$ interior collocation points, and $8000$ points for ring sampling near the boundaries. The corresponding performance metrics are reported in Table~\ref{FP_1+1}. The single-network PINN already achieves a mean MAE of $1.07\times 10^{-2}$, with the best run (seed $48$) reaching MAE $6.6\times 10^{-3}$.

\begin{table}[H]
\centering
\caption{One-net, one-phase PINN with ALM for the 1D Fokker--Planck equation. Performance metrics (mean $\pm$ standard deviation) across different seed values.}
\label{FP_1+1}
\resizebox{\textwidth}{!}{%
\begin{tabular}{@{}lccccccc@{}}
\toprule
 & \multicolumn{4}{c}{Interior} & \multicolumn{1}{c}{Boundary} & \multicolumn{1}{c}{PDE resid}  \\
\cmidrule(lr){2-5} \cmidrule(lr){6-6} \cmidrule(lr){7-7}
     & MAE & RMSE & Rel. $L^2$ & Accuracy (\%) & $L^2$ & $L^2(\text{residual})$ \\
\midrule
Seed $42$ & $0.0129$ & $0.0141$ & $0.0485$ & $95.15$ & $0.0119$ & $0.0490$ \\
Seed $44$  & $0.0114$ & $0.0153$  & $0.0526$ & $94.74$ & $0.0106$ & $0.0405$ \\
Seed $46$  & $0.0111$ & $0.0131$  & $0.0449$ & $95.51$ & $0.0146$ & $0.0542$ \\
Seed $48$ & $0.0066$ & $0.0078$ & $0.0267$ & $97.33$ & $0.0123$ & $0.0461$ \\
Seed $50$ & $0.0116$ & $0.0138$ & $0.0475$ & $95.25$ & $0.0151$ & $0.0626$ \\
\textbf{Mean $\pm$ Std} & $0.01072 \pm 0.00215$ & $0.01282 \pm 0.00261$ & $0.04404 \pm 0.00902$ & $95.60 \pm 0.90$ & $0.01290 \pm 0.00170$ & $0.05048 \pm 0.00751$
\\
\bottomrule
\end{tabular}
}
\end{table}

\paragraph{One-net and two-phase training.}
We next consider one-net and two-phase training within the same ALM-based PINN architecture. The network architecture is unchanged (six hidden layers with $50$ neurons). We train $10{,}000$ epochs in each of the two phases, again with early stopping. We use $500$ boundary points, $7000$ collocation points for the first phase, $8000$ for the second phase, and $8000$ points for ring sampling. The results in Table~\ref{FP_1+2} show that this curriculum slightly increases the mean error on this relatively simple 1D benchmark, and does not clearly outperform the one-phase baseline.

\begin{table}[H]
\centering
\caption{One-net, two-phase PINN with ALM for the 1D Fokker--Planck equation. Performance metrics (mean $\pm$ standard deviation) across different seed values.}
\label{FP_1+2}
\resizebox{\textwidth}{!}{%
\begin{tabular}{@{}lccccccc@{}}
\toprule
 & \multicolumn{4}{c}{Interior} & \multicolumn{1}{c}{Boundary} & \multicolumn{1}{c}{PDE resid}  \\
\cmidrule(lr){2-5} \cmidrule(lr){6-6} \cmidrule(lr){7-7}
     & MAE & RMSE & Rel. $L^2$ & Accuracy (\%) & $L^2$ & $L^2(\text{residual})$   \\
\midrule
Seed $42$ & $0.0143$ & $0.0181$ & $0.0622$ & $93.78$ & $0.0086$ & $0.0434$ \\
Seed $44$  & $0.0165$ & $0.0186$  & $0.0639$ & $93.61$ & $0.0221$ & $0.0719$ \\
Seed $46$  & $0.0164$ & $0.0227$  & $0.0779$ & $92.21$ & $0.0056$ & $0.0469$ \\
Seed $48$ & $0.0067$ & $0.0094$ & $0.0321$ & $96.79$ & $0.0075$ & $0.0309$ \\
Seed $50$ & $0.0131$ & $0.0156$ & $0.0537$ & $94.63$ & $0.0095$ & $0.0397$ \\
\textbf{Mean $\pm$ Std} & $0.01340 \pm 0.00359$ & $0.01688 \pm 0.00438$ & $0.05796 \pm 0.04670$ & $94.20 \pm 1.51$ & $0.01066 \pm 0.00587$ & $0.04656 \pm 0.01373$
\\
\bottomrule
\end{tabular}
}
\end{table}

\paragraph{Two-nets and two-phase training (original ALM setting).}
We then investigate a dual-network PINN with two-phase training under the original ALM formulation. The larger network (domain network) has one input layer, six hidden layers with $50$ neurons each, and one output layer, whereas the smaller network (boundary network) has one input layer, three hidden layers with $25$ neurons each, and one output layer. For each phase, we train up to $10{,}000$ epochs with early stopping. We use $500$ boundary points, $7000$ collocation points in the first phase, $8000$ in the second phase, and $8000$ points for ring sampling. The results in Table~\ref{FP_2+2} show that, for this 1D test, the dual-network ALM configuration does not improve over the simpler one-net baseline; in fact, the mean MAE slightly increases. This suggests that, in the absence of additional guidance, the extra flexibility and regularization introduced by the dual architecture can make the optimization problem harder rather than easier.

\begin{table}[H]
\centering
\caption{Two-nets, two-phase PINN with ALM for the 1D Fokker--Planck equation. Performance metrics (mean $\pm$ standard deviation) across different seed values.}
\label{FP_2+2}
\resizebox{\textwidth}{!}{%
\begin{tabular}{@{}lccccccc@{}}
\toprule
 & \multicolumn{4}{c}{Interior} & \multicolumn{1}{c}{Boundary} & \multicolumn{1}{c}{PDE resid}  \\
\cmidrule(lr){2-5} \cmidrule(lr){6-6} \cmidrule(lr){7-7}
     & MAE & RMSE & Rel. $L^2$ & Accuracy (\%) & $L^2$ & $L^2(\text{residual})$ \\
\midrule
Seed $42$ & $0.0184$ & $0.0244$ & $0.0838$ & $91.62$ & $0.0085$ & $0.0466$ \\
Seed $44$ & $0.0162$ & $0.0196$  & $0.0673$ & $93.27$ & $0.0159$ & $0.0768$ \\
Seed $46$ & $0.0136$ & $0.0182$  & $0.0624$ & $93.76$ & $0.0065$ & $0.0589$ \\
Seed $48$ & $0.0198$ & $0.0224$ & $0.0771$ & $92.29$ & $0.0206$ & $0.0674$ \\
Seed $50$ & $0.0105$ & $0.0119$ & $0.0407$ & $95.93$ & $0.0110$ & $0.0594$ \\
\textbf{Mean $\pm$ Std} & $0.0157 \pm 0.0036$ & $0.0193 \pm 0.0044$ & $0.0663 \pm 0.0160$ & $93.37 \pm 1.57$ & $0.0125 \pm 0.0055$ & $0.0618 \pm 0.0109$
\\
\bottomrule
\end{tabular}
}
\end{table}

\begin{figure}[H]
    \centering
    \includegraphics[width=\linewidth]{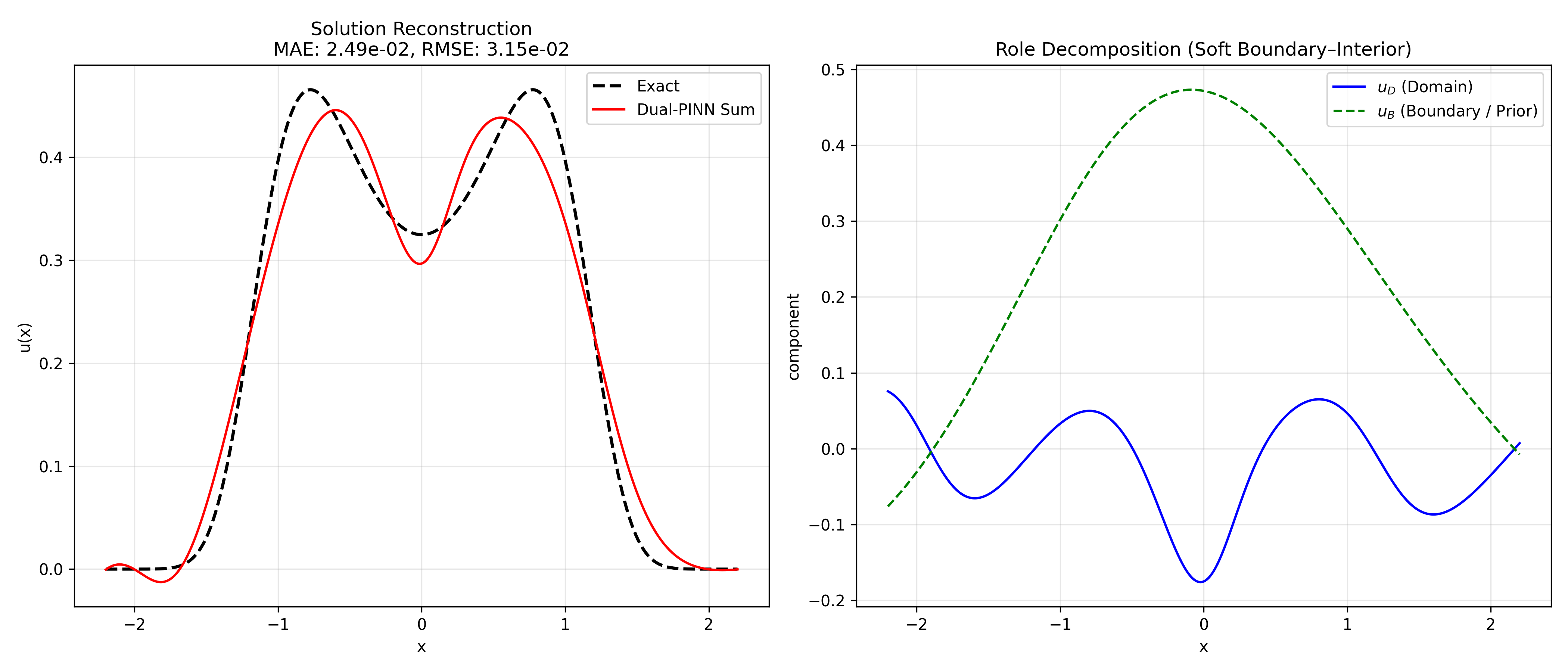}
    \caption{Dual-PINN with two networks and soft boundary--interior roles for the 1D Fokker--Planck equation (configuration of Table~\ref{FP_2+2}). Left: comparison between the exact solution (black dashed) and the combined prediction $u_{\mathrm D}+u_{\mathrm B}$ (red). Right: role decomposition. The ``boundary'' subnetwork $u_{\mathrm B}$ (green) remains active across the whole domain, while the ``domain'' subnetwork $u_{\mathrm D}$ (blue) produces small oscillatory corrections, showing that the desired specialization does not emerge in this joint ALM training.}
    \label{fig:FP-dual-soft-failure}
\end{figure}

Figure~\ref{fig:FP-dual-soft-failure} visualizes a representative run from this two-net ALM setting. The boundary subnetwork $u_{\mathrm B}$ ends up modeling most of the stationary density, while the domain subnetwork $u_{\mathrm D}$ only contributes small oscillatory corrections. Soft boundary-interior specialization fails to emerge when both subnetworks are trained jointly under competing losses, which explains the lack of improvement seen in Table~\ref{FP_2+2}.

\paragraph{Dual-network residual PINN with sequential training.}
Motivated by these observations, we further explore a variant of the dual-network idea tailored to the Fokker--Planck benchmark. Instead of training both subnetworks simultaneously under competing losses, we decompose the solution as
\begin{equation}
u(x) = u_{\mathrm B}(x) + u_{\mathrm D}(x),
\end{equation}
where $u_{\mathrm B}$ plays the role of a data-driven baseline and $u_{\mathrm D}$ acts as a physics-based residual. Training proceeds in three sequential phases:

\begin{itemize}
    \item \textbf{Phase~1 (data prior).} We pretrain $u_{\mathrm B}$ on pseudo-measurements sampled from the closed-form solution, minimizing a mean squared error loss on pairs $(x_i,u(x_i))$ distributed over $[a,b]$. The domain network $u_{\mathrm {D}} $ remains inactive in this phase. This yields a smooth baseline that already approximates the stationary density well.
    \item \textbf{Phase~2 (physics residual).} The parameters of $u_{\mathrm B}$ are frozen, and $u_{\mathrm D}$ is optimized to reduce the Fokker--Planck residual and the normalization error of the combined prediction $u_{\mathrm B}+u_{\mathrm D}$. We also add a mild pinning term that keeps $u_{\mathrm D}$ small at the domain boundaries, so that the near-boundary behavior remains dominated by the baseline.
    \item \textbf{Phase~3 (joint fine-tuning).} Finally, we jointly update both subnetworks using an objective that combines the PDE residual, the normalization constraint, and the data loss on the pseudo-measurements. In this phase, we also include a soft role prior that penalizes the energy of $u_{\mathrm B}$ in the interior and of $u_{\mathrm D}$ near the boundary, encouraging a boundary-interior specialization without enforcing hard masks.
\end{itemize}

For this dual-residual PINN, we reuse the same backbone architectures as in the two-net ALM setting (six hidden layers with $50$ neurons for $u_{\mathrm D}$, and a smaller network for $u_{\mathrm B}$). A single run of the sequential dual-residual scheme yields an interior MAE of approximately $2.8\times 10^{-3}$ and an RMSE of $3.2\times 10^{-3}$ on a fine evaluation grid, which is about $2.4\times$ lower MAE than the best one-net baseline (seed $48$) and almost $4\times$ lower than the mean MAE reported in Table~\ref{FP_1+1}. The comparison is summarized in Table~\ref{FP_dual_residual}.

\begin{table}[H]
\centering
\caption{Comparison between the best one-net PINN and the dual-network residual PINN (sequential training) on the 1D Fokker--Planck equation.}
\label{FP_dual_residual}
\begin{tabular}{@{}lcc@{}}
\toprule
Method & MAE (interior) & RMSE (interior) \\
\midrule
One-net, one-phase (best seed $48$) & $6.6\times 10^{-3}$ & $7.8\times 10^{-3}$ \\
Dual-network residual (sequential)  & $2.8\times 10^{-3}$ & $3.2\times 10^{-3}$ \\
\bottomrule
\end{tabular}
\end{table}

\begin{figure}[H]
    \centering
    \includegraphics[width=0.7\linewidth]{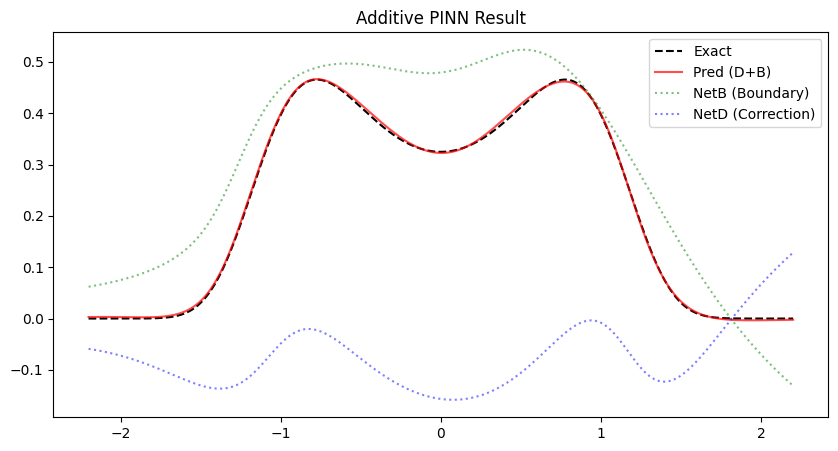}
    \caption{Dual-network residual PINN for the 1D Fokker--Planck equation. The black dashed line is the analytical solution, while the red line shows the combined prediction $u_{\mathrm B}+u_{\mathrm D}$. The green and blue curves indicate respectively the contributions of the baseline subnetwork $u_{\mathrm B}$ and of the residual subnetwork $u_{\mathrm D}$. The residual network primarily refines the internal structure, while the baseline controls the overall shape.}
    \label{fig:FP-dual-residual}
\end{figure}

Figure~\ref{fig:FP-dual-residual} shows the prediction of the dual-network residual PINN compared to the analytical solution, together with the contributions of $u_{\mathrm B}$ and $u_{\mathrm D}$. The baseline subnetwork captures the smooth global profile, while the residual subnetwork contributes small corrections that sharpen the double-peaked structure and enforce the PDE and normalization constraints. This experiment demonstrates that decomposing the solution into a data-driven prior and a physics-based residual, and training the two components sequentially, can substantially reduce the error, even on a relatively benign 1D benchmark where naive dual-network ALM does not provide significant benefits.

\subsection{2D Laplace Equation}\label{6.1}
\begin{figure}[H]
    \centering
     \begin{subfigure}[b]{\textwidth}
         \centering
         \includegraphics[height=4.2cm,width=\textwidth]{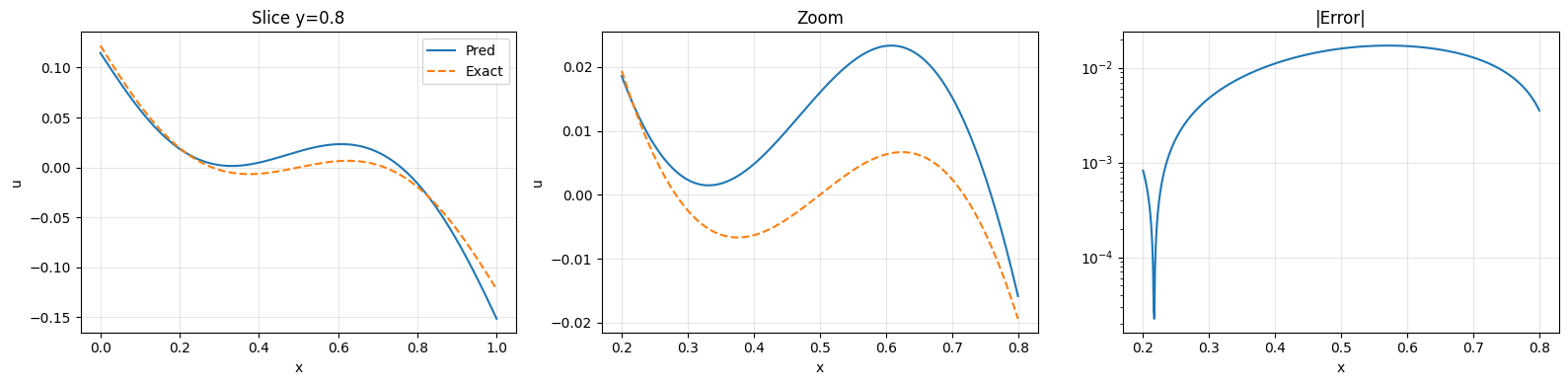}
         \caption{Predicted vs. exact solution at $y=0.8$ with zoomed view and error.}
     \end{subfigure}
     ~~
     \begin{subfigure}[b]{\textwidth}
         \centering
         \includegraphics[height=4.2cm,width=\textwidth]{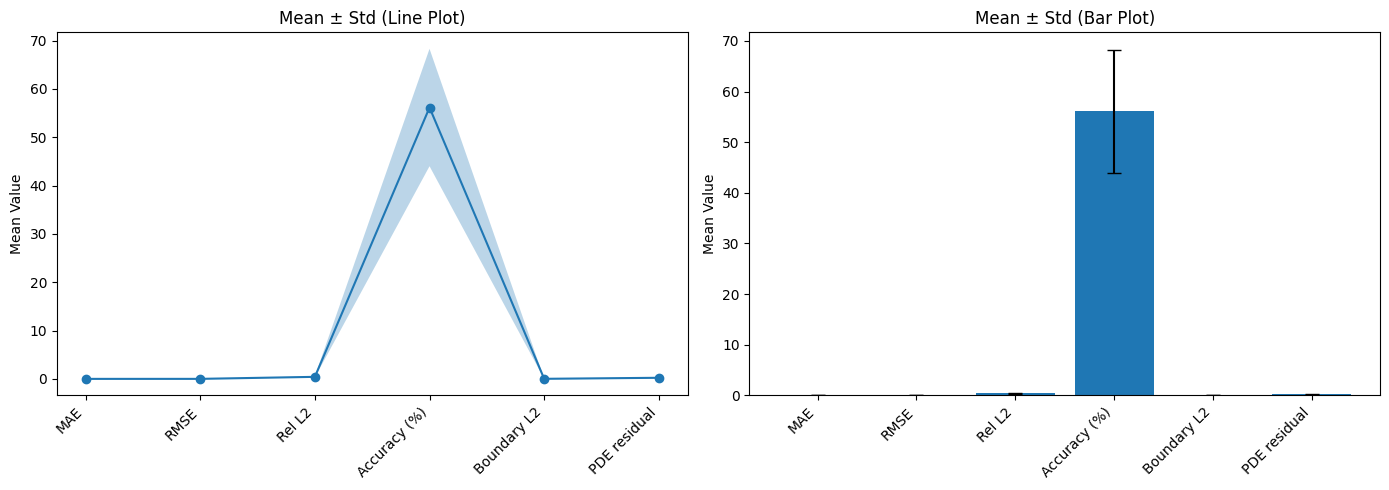}
         \caption{Mean $\pm$ std of evaluation metrics (line and bar plots).}
     \end{subfigure}
    \caption{Comparison along slice $y=0.8$ and statistical summary of evaluation metrics with seed $46$.}\label{LE_1+1}
\end{figure}

In this subsection, we solve the 2D \textbf{Laplace equation} on the unit square $[0,1]^2$:
\begin{equation}
\nabla^2 u(x,y) \equiv u_{xx}(x,y) + u_{yy}(x,y) = 0,
\quad \text{for}~~(x,y) \in (0,1)^2
\end{equation}
with the non-trivial boundary conditions $u(x,y) = u_{\text{exact}}(x,y)$ \text{on the boundary } $\partial \Omega$. There is no source term (forcing), i.e.
$f=0$. Then the exact solution is given by
\begin{equation}
u_{\text{exact}}(x,y)
= \frac{1}{16}\left( x^5 - 10 x^3 y^2 + 5 xy^4 \right).
\end{equation}
This polynomial is harmonic (its Laplacian is zero), so it is a valid exact solution to the Laplace equation. The boundary values are simply this exact function restricted to $\partial \Omega$. 
\paragraph{One-net and one-phase training.}
We begin by solving the Laplace equation on the domain $[0,1]\times[0,1]$
using a standard PINNs framework. Specifically, we adopt a single-network, single-phase training strategy, with boundary conditions imposed through the augmented Lagrangian method (ALM). The neural network architecture consists of seven hidden layers with $64$ neurons each, $\tanh$ activation, and Xavier initialization. Boundary points are generated along each of the four edges of the square using Latin hypercube sampling (LHS), with 300 points per edge. Additionally, 7000 collocation points are sampled uniformly within the interior of the domain. Training proceeds for up to 2000 epochs, using early stopping to promote stability, reduce computational cost, and mitigate overfitting. Table \eqref{tab:onenet1} reports the evaluation metrics including MAE, RMSE, relative $L^2$ error, accuracy, BC $L^2$ and PDE $L^2$. The single-network PINN achieves a mean MAE of $1.7\times 10^{-2}$, with the best run (seed 46) reaching $9.9\times 10^{-3}$. Figure \eqref{LE_1+1} presents the solution comparison at $y=0.8$, along with the corresponding statistical metrics.

\begin{table}[H]
\centering
\caption{One-net, one-phase PINN with ALM for Laplace equation. Performance metrics (mean $\pm$ standard deviation) across different seed values.}
\label{tab:onenet1}
\resizebox{\textwidth}{!}{%
\begin{tabular}{@{}lccccccccl@{}}
\toprule
 & \multicolumn{4}{c}{Interior} & \multicolumn{1}{c}{Boundary} & \multicolumn{1}{c}{PDE resid}  \\
\cmidrule(lr){2-5} \cmidrule(lr){6-6} \cmidrule(lr){7-7}
     & MAE & RMSE & Rel. $L^2$ & Accuracy (\%) & $L^2$ & $L^2$(residual) &   \\
\midrule
Seed $40$ & $0.0218$ & $0.0239$ & $0.4958$ & $50.42$ & $0.0261$ & $0.2534$   \\
Seed $42$  & $0.0214$ & $0.0258$  & $0.5351$ & $46.49$ & $0.0397$ & $0.2681$ \\
Seed $44$  & $0.0138$ & $0.0174$  & $0.3612$ & $63.88$ & $0.0299$ & $0.3911$ \\
Seed $46$ & $0.0099$ & $0.0127$ & $0.2636$ & $73.64$ & $0.0226$ & $0.1565$
\\
Seed $48$ & $0.0223$ & $0.0259$ & $0.5372$ & $46.28$ & $0.0370$ & $0.1554$ 
\\
\textbf{Mean $\pm$ Std} & $0.0178 \pm 0.0056$ & $0.0211 \pm 0.0059$ & $0.4386 \pm 0.1214$ & $56.14 \pm 12.14$ & $0.0311 \pm 0.0072$ & $0.2449 \pm 0.0972$ \\
\bottomrule
\end{tabular}
}
\end{table}


\begin{figure}[H]
    \centering
     \begin{subfigure}[b]{\textwidth}
         \centering
         \includegraphics[height=4.2cm,width=\textwidth]{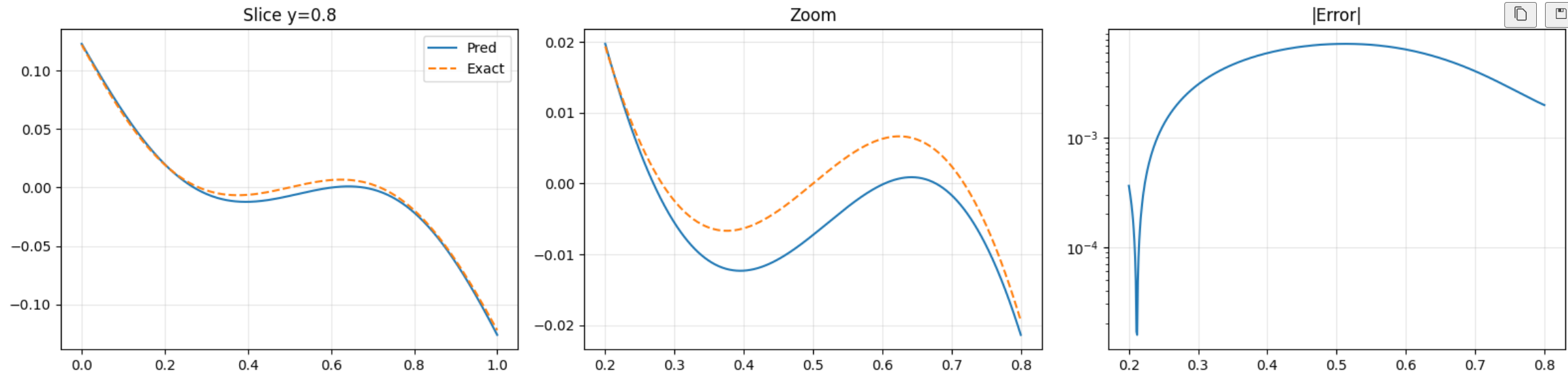}
         \caption{Predicted vs. exact solution at $y=0.8$ with zoomed view and error.}
     \end{subfigure}
     ~~
     \begin{subfigure}[b]{\textwidth}
         \centering
         \includegraphics[height=4.2cm,width=\textwidth]{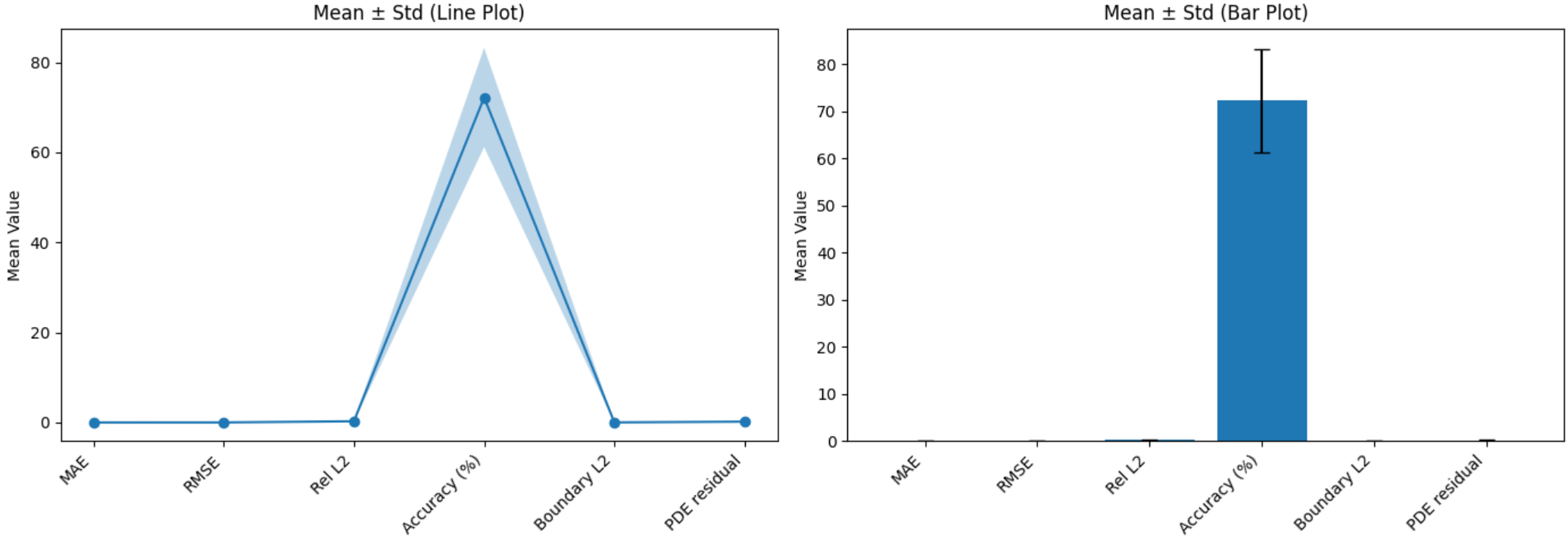}
         \caption{Mean $\pm$ std of evaluation metrics (line and bar plots).}
     \end{subfigure}
    \caption{Comparison along slice $y=0.8$ and statistical summary of evaluation metrics with seed $48$.}\label{LE_1+2}
\end{figure}

\paragraph{One-net and two-phase training.}
Next, we employ the same PINNs framework based on ALM by implementing the \enquote{one-net and two-phase training}. The network architecture is the same (seven hidden layers, each with 64 neurons). Boundary points are generated along each of the four edges of the square using LHS. Also, a \enquote{ring sampling}, denser near the boundary, is added to focus on near-boundary accuracy. We train $1000$ epochs in each of two phases with early stopping. We generate $ 7,000$ uniform collocation points in the first phase and $ 8,000$ collocation points with ring sampling in the second phase. The results in the Table \eqref{tab:onenet} show that this setup increases MAE (mean MAE $1.0\times 10^{-2}$ and MAE $4.8\times 10^{-3}$) for the considered equation. The comparison at $y=0.8$ and the associated statistical metrics are shown in Figure \eqref{LE_1+2}.

\begin{table}[H]
\centering
\caption{One-net, two-phase PINN with ALM for Laplace equation. Performance metrics (mean $\pm$ standard deviation) across different seed values.}\label{tab:onenet}
\resizebox{\textwidth}{!}{%
\begin{tabular}{@{}lccccccccl@{}}
\toprule
 & \multicolumn{4}{c}{Interior} & \multicolumn{1}{c}{Boundary} & \multicolumn{1}{c}{PDE resid}  \\
\cmidrule(lr){2-5} \cmidrule(lr){6-6} \cmidrule(lr){7-7}
     & MAE & RMSE & Rel. $L^2$ & Accuracy (\%) & $L^2$ & $L^2$(residual) &   \\
\midrule
Seed $40$ & $0.0101$ & $0.0131$ & $0.2707$ & $72.93$ & $0.0243$ & $0.3099$   \\
Seed $42$  & $0.0146$ & $0.0216$  & $0.4481$ & $55.19$ & $0.0418$ & $0.1030$  \\
Seed $44$  & $0.0092$ & $0.0112$  & $0.2319$ & $76.81$ & $0.0139$ & $0.2547$ \\
Seed $46$ & $0.0123$ & $0.0141$ & $0.2922$ & $70.78$ & $0.0176$ & $0.1390$\\
Seed $48$ & $0.0048$ & $0.0071$ & $0.1464$ & $85.36$ & $0.0151$ & $0.1425$ \\
\textbf{Mean $\pm$ Std} & $0.0102 \pm 0.0037$ & $0.0134 \pm 0.0053$ & $0.2779 \pm 0.1103$ & $72.21 \pm 11.03$ & $0.0225 \pm 0.0115$ & $0.1898 \pm 0.0880$ \\
\bottomrule
\end{tabular}
}
\end{table}

\begin{figure}[H]
    \centering
     \begin{subfigure}[b]{\textwidth}
         \centering
         \includegraphics[height=4.2cm,width=\textwidth]{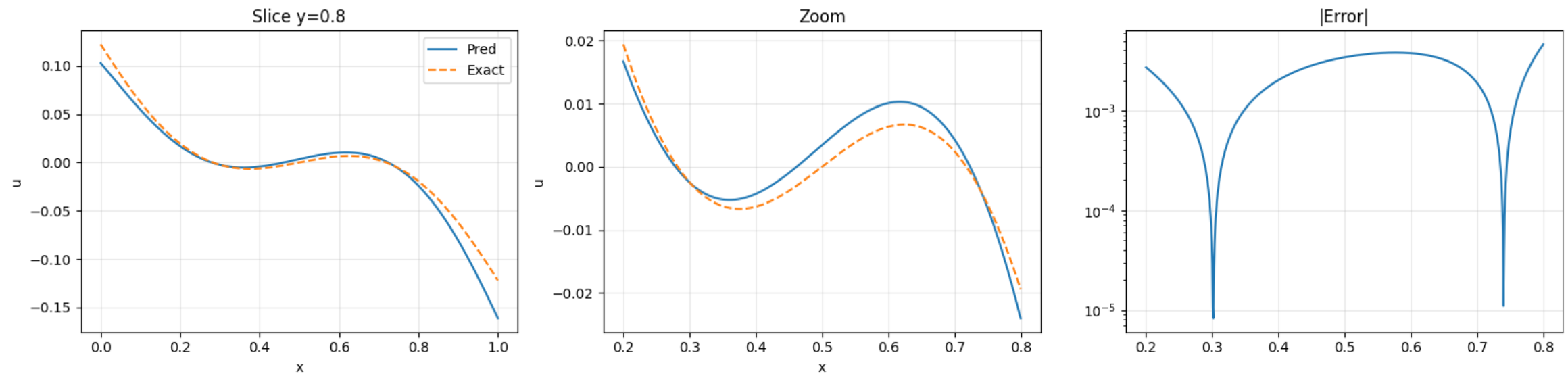}
         \caption{Predicted vs. exact solution at $y=0.8$ with zoomed view and error.}
     \end{subfigure}
     ~~
     \begin{subfigure}[b]{\textwidth}
         \centering
         \includegraphics[height=4.2cm,width=\textwidth]{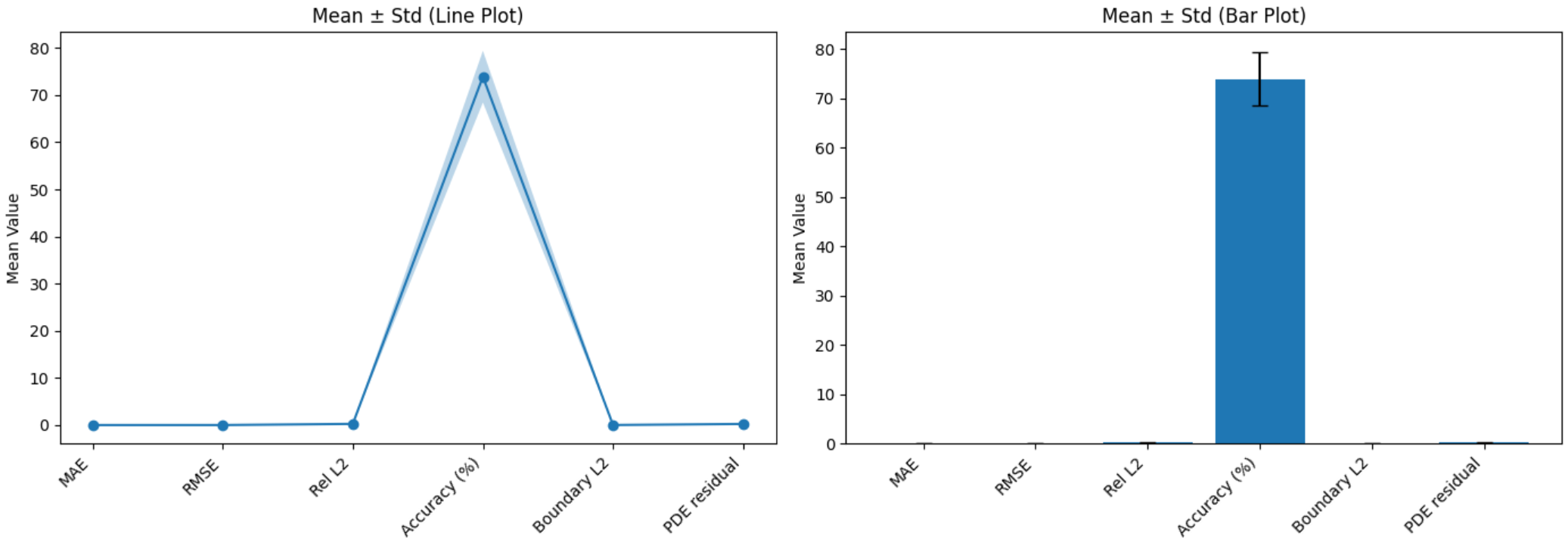}
         \caption{Mean $\pm$ std of evaluation metrics (line and bar plots).}
     \end{subfigure}
    \caption{Comparison along slice $y=0.8$ and statistical summary of evaluation metrics with seed $42$.}\label{LE_2+2}
\end{figure}

\paragraph{Two-nets and two-phase training.}
We now examine a dual-network PINNs framework based on the ALM formulation. The proposed architecture consists of two neural networks: a domain network (larger) with six hidden layers of 48 neurons each, and a boundary network (smaller) with three hidden layers of 16 neurons each. We generate 1000 boundary points, 6000 collocation points, and 8000 ring-sampling points. Training is carried out for 2000 epochs in each of the two phases, with early stopping applied in both.
In Phase 1, we use 7000 uniformly distributed collocation points together with strong ALM-enforced boundary conditions. In Phase 2, we employ 8000 collocation points following a ring-sampling strategy combined with annealing.
Table \eqref{tab:twonet} summarizes the performance metrics, where we obtain a mean MAE of $9.6\times 10^{-3}$ and and a best-run MAE of $4.5\times 10^{-3}$. Figure \eqref{LE_2+2} shows the solution comparison along $y=0.8$, together with the corresponding statistical indicators.

\begin{table}[H]
\centering
\caption{Two-nets, two-phase PINN with ALM for Laplace equation. Performance metrics (mean $\pm$ standard deviation) across different seed values.}
\label{tab:twonet}
\resizebox{\textwidth}{!}{%
\begin{tabular}{@{}lccccccccl@{}}
\toprule
 & \multicolumn{4}{c}{Interior} & \multicolumn{1}{c}{Boundary} & \multicolumn{1}{c}{PDE resid}  \\
\cmidrule(lr){2-5} \cmidrule(lr){6-6} \cmidrule(lr){7-7}
     & MAE & RMSE & Rel. $L^2$ & Accuracy (\%) & $L^2$ & $L^2$(residual) &   \\
\midrule
Seed $40$ & $0.0107$ & $0.0129$ & $0.2651$ & $73.49$ & $0.0197$ & $0.2031$  \\
Seed $42$  & $0.0045$ & $0.0082$  & $0.1691$ & $83.09$ & $0.0176$ & $0.3151$
\\
Seed $44$ & $0.0096$ & $0.0128$  & $0.2636$ & $73.64$ & $0.0184$ & $0.2061$ \\
Seed $46$ & $0.0131$ & $0.0148$ & $0.3040$ & $69.60$ & $0.0169$ & $0.2360$\\
Seed $48$ & $0.0104$ & $0.0147$ & $0.3027$ & $69.73$ & $0.0318$ & $0.2759$ \\
\textbf{Mean $\pm$ Std} & $0.00966 \pm 0.00317$ & $0.01268 \pm 0.00268$ & $0.26090 \pm 0.05490$ & $73.91 \pm 5.49$ & $0.02088 \pm 0.00619$ & $0.24724 \pm 0.04794$ \\
\bottomrule
\end{tabular}
}
\end{table}

The MAE, RMSE, and mean values for each approach on the Laplace equation are compared in the following Table \eqref{comparison1}. These findings verify that the dual-network ALM setup outperforms the more straightforward one-net with one and two phase training baselines for this 2D test example.

\begin{table}[H]
\centering
\caption{Comparison between the one-net PINN with one and two phases and the Dual-network PINN on the Laplace equation.}
\label{comparison1}
\begin{tabular}{@{}lccc@{}}
\toprule
Method & MAE (interior) & RMSE (interior) & Mean (MAE) \\
\midrule
One-net, one-phase (best seed $46$) & $9.9\times 10^{-3}$ & $1.2\times 10^{-2}$ & $1.7\times 10^{-2}$ \\
One-net, two-phase (best seed $48$) & $4.8\times 10^{-3}$ & $7.1\times 10^{-3}$ & $1.0\times 10^{-2}$\\
Dual-network, two-phase (best seed $42$)   & $4.5\times 10^{-3}$ & $8.2\times 10^{-3}$ & $9.6\times 10^{-3}$\\
\bottomrule
\end{tabular}
\end{table}

The performance metrics in Table \eqref{tab:Role-Priors}, obtained using the Dual-PINN without role priors, show a substantial degradation across all measures. In comparison with our Dual-PINN with role priors (Table \ref{tab:twonet}), the results are significantly worse for every metric.

\begin{table}[H]
\centering
\caption{Two-nets, two-phase PINN without role priors with ALM for Laplace equation. Performance metrics across different seed values.}
\label{tab:Role-Priors}
\resizebox{\textwidth}{!}{%
\begin{tabular}{@{}lccccccccl@{}}
\toprule
 & \multicolumn{4}{c}{Interior} & \multicolumn{1}{c}{Boundary} & \multicolumn{1}{c}{PDE resid}  \\
\cmidrule(lr){2-5} \cmidrule(lr){6-6} \cmidrule(lr){7-7}
     & MAE & RMSE & Rel. $L^2$ & Accuracy (\%) & $L^2$ & $L^2$(residual) &   \\
\midrule
Seed $40$ & $0.0185$ & $0.0227$ & $0.4671$ & $53.29$ & $0.0375$ & $0.3186$  \\
Seed $42$  & $0.0116$ & $0.0144$  & $0.2958$ & $70.42$ & $0.0215$ & $0.2680$
\\
Seed $44$ & $0.0110$ & $0.0148$  & $0.3058$ & $69.42$ & $0.0240$ & $0.2352$ \\
Seed $46$ & $0.0094$ & $0.0147$ & $0.3033$ & $69.67$ & $0.0223$ & $0.2370$\\
Seed $48$ & $0.0088$ & $0.0118$ & $0.2422$ & $75.78$ & $0.0234$ & $0.2743$ \\
\textbf{Mean $\pm$ Std} & $0.01186 \pm 0.00388$ & $0.01568 \pm 0.00411$ & $0.32284 \pm 0.08473$ & $67.72 \pm 8.47$ & $0.02574 \pm 0.00664$ & $0.26662 \pm 0.03401$ \\
\bottomrule
\end{tabular}
}
\end{table}


\subsection{2D Poisson Equation}
The equation is commonly encountered in fluid dynamics.The considering equation is not varying in time, but still, hard to solve analytically and therefore a numerical approach is usually required.
\begin{equation}
u_{xx} + u_{yy} = f(x,y), \quad x \in [0,1], \quad y \in [0,1]
\end{equation}
with the non-trivial boundary conditions $u(x,y) = u_{\text{exact}}(x,y)$ \text{on the boundary} $\partial \Omega$. There is a source term (forcing), i.e. $f(x,y)=-\sin(\pi x) \sin(\pi y)$. Then the exact solution is given by
\begin{equation}
u_{\text{exact}}(x, y) = \frac{1}{2\pi^{2}} \sin(\pi x)\,\sin(\pi y).
\end{equation}

\paragraph{One-net and one-phase training.}
We apply the standard PINNs framework to solve the Poisson equation using a single-network, single-phase training strategy with ALM-enforced boundary conditions. Boundary points are generated using LHS, with 300 points on each edge of the square domain. Additionally, 7000 collocation points are sampled uniformly throughout the interior. The neural network architecture remains the same as in the Laplace case-seven hidden layers with 64 neurons each. Training proceeds for up to 2000 epochs, with early stopping applied to stabilize learning and reduce computational cost. Table \eqref{PE_1+1} reports the performance metrics, where the method achieves a mean MAE of $4.5\times 10^{-3}$ and a best-run MAE of $2.8\times 10^{-3}$ (seed 46). 

\begin{table}[H]
\centering
\caption{One-net, one-phase PINN with ALM for Poisson equation. Performance metrics (mean $\pm$ standard deviation) across different seed values.}
\label{PE_1+1}
\resizebox{\textwidth}{!}{%
\begin{tabular}{@{}lccccccccl@{}}
\toprule
 & \multicolumn{4}{c}{Interior} & \multicolumn{1}{c}{Boundary} & \multicolumn{1}{c}{PDE resid}  \\
\cmidrule(lr){2-5} \cmidrule(lr){6-6} \cmidrule(lr){7-7}
     & MAE & RMSE & Rel. $L^2$ & Accuracy (\%) & $L^2$ & $L^2$(residual) & \\
\midrule
Seed $40$ & $0.0046$ & $0.0053$ & $0.2122$ & $78.78$ & $0.0079$ & $0.5042$ \\
Seed $42$ & $0.0037$ & $0.0046$  & $0.1837$ & $81.63$ & $0.0056$ & $0.4982$ \\
Seed $44$ & $0.0048$ & $0.0059$  & $0.2324$ & $76.76$ & $0.0062$ & $0.5226$ \\
Seed $46$ & $0.0028$ & $0.0035$ & $0.1397$ & $86.03$ & $0.0041$ & $0.4958$ \\
Seed $48$ & $0.0069$ & $0.0083$ & $0.3287$ & $67.13$ & $0.0091$ & $0.4918$ \\
\textbf{Mean $\pm$ Std} & $0.00456 \pm 0.00153$ & $0.00552 \pm 0.00179$ & $0.21934 \pm 0.07035$ & $78.07 \pm 7.03$ & $0.00658 \pm 0.00196$ & $0.50252 \pm 0.01209$
\\
\bottomrule
\end{tabular}
}
\end{table}

\paragraph{One-net and two-phase training.}
Next, we apply the same ALM-based PINNs framework to solve the Poisson equation using a single-network, two-phase training strategy. As before, we generate 300 boundary points using LHS and employ ring sampling to increase point density near the boundary and improve boundary accuracy. The network architecture remains unchanged, consisting of seven hidden layers with 64 neurons each. Training is performed for 1000 epochs in each phase, with early stopping activated to enhance stability. In Phase 1, we use 7000 uniformly distributed collocation points, while in Phase 2 we switch to 8000 collocation points generated through ring sampling.
Table \eqref{PE_1+2} summarizes the performance, showing a mean MAE of
$2.4\times 10^{-3}$) and a best-run MAE of $1.1\times 10^{-3}$ (seed $48$).

\begin{table}[H]
\centering
\caption{One-net, two-phase PINN with ALM for Poisson equation. Performance metrics (mean $\pm$ standard deviation) across different seed values.}
\label{PE_1+2}
\resizebox{\textwidth}{!}{%
\begin{tabular}{@{}lccccccccl@{}}
\toprule
 & \multicolumn{4}{c}{Interior} & \multicolumn{1}{c}{Boundary} & \multicolumn{1}{c}{PDE resid}  \\
\cmidrule(lr){2-5} \cmidrule(lr){6-6} \cmidrule(lr){7-7}
     & MAE & RMSE & Rel. $L^2$ & Accuracy (\%) & $L^2$ & $L^2$(residual) & \\
\midrule
Seed $40$ & $0.0033$ & $0.0040$ & $0.1602$ & $83.98$ & $0.0058$ & $0.4961$ \\
Seed $42$ & $0.0011$ & $0.0012$  & $0.0494$ & $95.06$ & $0.0015$ & $0.4987$ \\
Seed $44$ & $0.0018$ & $0.0022$  & $0.0885$ & $91.15$ & $0.0035$ & $0.5055$ \\
Seed $46$ & $0.0050$ & $0.0063$ & $0.2517$ & $74.83$ & $0.0086$ & $0.4940$ \\
Seed $48$ & $0.0011$ & $0.0014$ & $0.0572$ & $94.28$ & $0.0020$ & $0.4978$ \\
\textbf{Mean $\pm$ Std} & $0.00246 \pm 0.00168$ & $0.00302 \pm 0.00214$ & $0.1214 \pm 0.08496$ & $87.86 \pm 8.50$ & $0.00428 \pm 0.00294$ & $0.49842 \pm 0.00434$
\\
\bottomrule
\end{tabular}
}
\end{table}

\begin{figure}[H]
    \centering
     \begin{subfigure}[b]{\textwidth}
     \centering         \includegraphics[height=3.9cm,width=\textwidth]{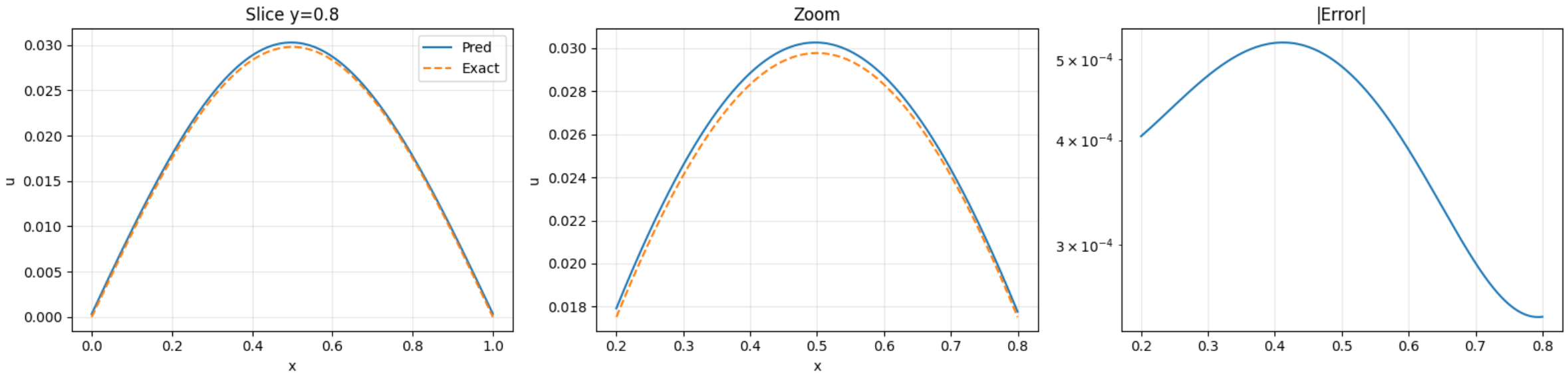}
     \caption{Predicted vs. exact solution at $y=0.8$ with zoomed view and error.}
     \end{subfigure}
     ~~
     \begin{subfigure}[b]{\textwidth}
    \centering
    \includegraphics[height=4.2cm,width=\textwidth]{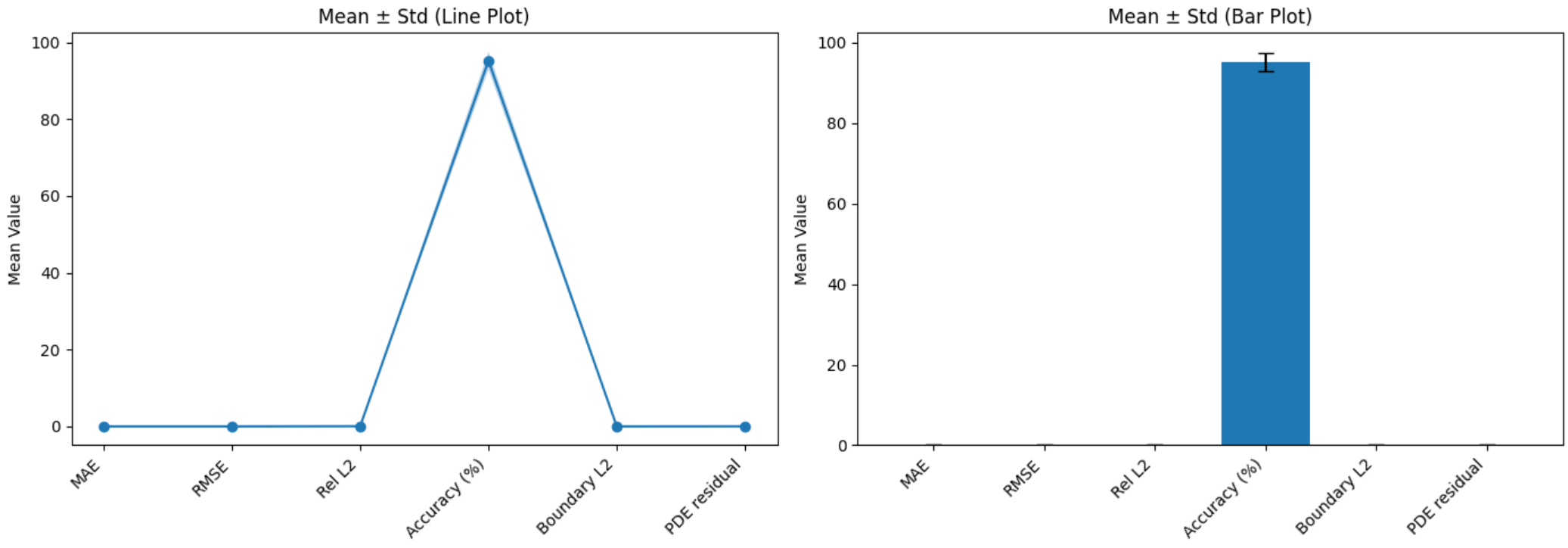}
    \caption{Mean $\pm$ std of evaluation metrics (line and bar plots).}
     \end{subfigure}
    \caption{Comparison along slice $y=0.8$ and statistical summary of evaluation metrics with seed $48$.}\label{FPE_2+2}
\end{figure}

\paragraph{Two-nets and two-phase training.}
We now examine our proposed framework: a dual-network PINNs approach with two-phase training based on the ALM formulation. The architecture consists of two neural networks: a larger domain network with six hidden layers of 48 neurons each, responsible for learning the interior solution, and a smaller boundary network with three hidden layers of 16 neurons each, designed to model the boundary conditions and capture near-boundary corrections. We employ 300 boundary points, along with 7000 collocation points in Phase 1 and 8000 collocation points in Phase 2. Additionally, 8000 ring-sampling points are used to densify samples near the boundary. Training runs for up to 2000 epochs in each phase, with early stopping criteria. Table \eqref{PE_2+2} reports the performance metrics, showing a mean MAE of $1.0\times 10^{-3}$) and a best-run MAE of $3.0\times 10^{-4}$ (seed $48$). Figure \eqref{FPE_2+2} presents the solution comparison along $y=0.8$.

\begin{table}[H]
\centering
\caption{Two-nets, two-phase PINN with ALM for Poisson equation. Performance metrics (mean $\pm$ standard deviation) across different seed values.}
\label{PE_2+2}
\resizebox{\textwidth}{!}{%
\begin{tabular}{@{}lccccccccl@{}}
\toprule
 & \multicolumn{4}{c}{Interior} & \multicolumn{1}{c}{Boundary} & \multicolumn{1}{c}{PDE resid}  \\
\cmidrule(lr){2-5} \cmidrule(lr){6-6} \cmidrule(lr){7-7}
     & MAE & RMSE & Rel. $L^2$ & Accuracy (\%) & $L^2$ & $L^2$(residual) & \\
\midrule
Seed $40$ & $0.0009$ & $0.0010$ & $0.0412$ & $95.88$ & $0.0010$ & $0.0123$ \\
Seed $42$ & $0.0008$ & $0.0010$  & $0.0387$ & $96.13$ & $0.0009$ & $0.0404$ \\
Seed $44$ & $0.0013$ & $0.0016$  & $0.0618$ & $93.82$ & $0.0021$ & $0.0148$ \\
Seed $46$ & $0.0018$ & $0.0021$ & $0.0842$ & $91.58$ & $0.0025$ & $0.0267$ \\
Seed $48$ & $0.0003$ & $0.0003$ & $0.0139$ & $98.61$ & $0.0005$ & $0.0115$ \\
\textbf{Mean $\pm$ Std} & $0.00102 \pm 0.00050$ & $0.00120 \pm 0.00061$ & $0.04796 \pm 0.02365$ & $95.20 \pm 2.36$ & $0.00140 \pm 0.00076$ & $0.02114 \pm 0.01107$
\\
\bottomrule
\end{tabular}
}
\end{table}

In the following Table \eqref{comparison}, we compare the MAE, RMSE and mean values for each method in the Poisson equation. These results confirm that for this 2D test, the dual-network ALM configuration improves over the simpler one-net with one- and two-phase training baselines.

\begin{table}[H]
\centering
\caption{Comparison between the one-net PINN with one and two phases and the dual-network PINN on the Poisson equation.}
\label{comparison}
\begin{tabular}{@{}lccc@{}}
\toprule
Method & MAE (interior) & RMSE (interior) & Mean (MAE) \\
\midrule
One-net, one-phase (best seed $46$) & $2.8\times 10^{-3}$ & $3.5\times 10^{-3}$ & $4.5\times 10^{-3}$ \\
One-net, two-phase (best seed $48$) & $1.1\times 10^{-3}$ & $1.4\times 10^{-3}$ & $2.4\times 10^{-3}$\\
Dual-network, two-phase (best seed $48$)   & $3.0\times 10^{-4}$ & $3.0\times 10^{-4}$ & $1.0\times 10^{-3}$\\
\bottomrule
\end{tabular}
\end{table}


\subsection{Wave equation}

In this final example, we consider the one-dimensional wave equation, which presents a significantly more challenging test case due to its multi-scale, time-dependent nature:
\begin{equation*}
\frac{\partial^2 u}{\partial t^2} - 4 \frac{\partial^2 u}{\partial x^2} = 0, \quad x \in [0, 1],\ t \in [0, 1],
\end{equation*}
with homogeneous Dirichlet boundary conditions
\begin{equation*}
u(0,t)=u(1,t)=0,\quad t\in[0, 1],
\end{equation*}
and initial conditions comprising two distinct Fourier modes:
\begin{equation*}
u(x, 0) = \sin(\pi x) + \frac{1}{2} \sin(4\pi x), \quad \frac{\partial u}{\partial t}(x, 0) = 0, \quad x \in [0, 1].
\end{equation*}
The exact solution is given by the superposition of two propagating modes:
\begin{equation}\label{eq:wave_exact}
u(x, t) = \sin(\pi x)\cos(2\pi t) + \frac{1}{2} \sin(4\pi x)\cos(8\pi t).
\end{equation}
This solution exhibits multi-scale behavior in both spatial ($\pi$ and $4\pi$) and temporal ($2\pi$ and $8\pi$) frequencies.

\textbf{Two-nets and two-phase training (baseline).} We first apply the base dual-network framework with tanh activations, following the same configuration used for the Laplace and Poisson benchmarks. The domain network consists of six hidden layers with 48 neurons, while the boundary network comprises three hidden layers with 16 neurons. Boundary and initial conditions are enforced via ALM, and training proceeds in two phases with uniform and ring sampling, respectively. The results in Table~\ref{tab:wave_baseline} show that the tanh-based configuration achieves only limited accuracy on this multi-scale benchmark, with MAE $1.88 \times 10^{-1}$ and accuracy $58.5\%$.

\begin{table}[H]
\centering
\caption{Two-nets, two-phase PINN with ALM for the wave equation using tanh activations (baseline).}
\label{tab:wave_baseline}
\begin{tabular}{lccccc}
\hline
& \multicolumn{3}{c}{Interior} & Boundary & \\
& MAE & RMSE & Rel.~$L^2$ & $L^2$ & Accuracy (\%) \\
\hline
Dual-net (tanh) & $1.88 \times 10^{-1}$ & $2.32 \times 10^{-1}$ & $4.15 \times 10^{-1}$ & $3.25 \times 10^{-2}$ & 58.5 \\
\hline
\end{tabular}
\end{table}

\begin{figure}[H]
    \centering
    \includegraphics[width=0.95\textwidth]{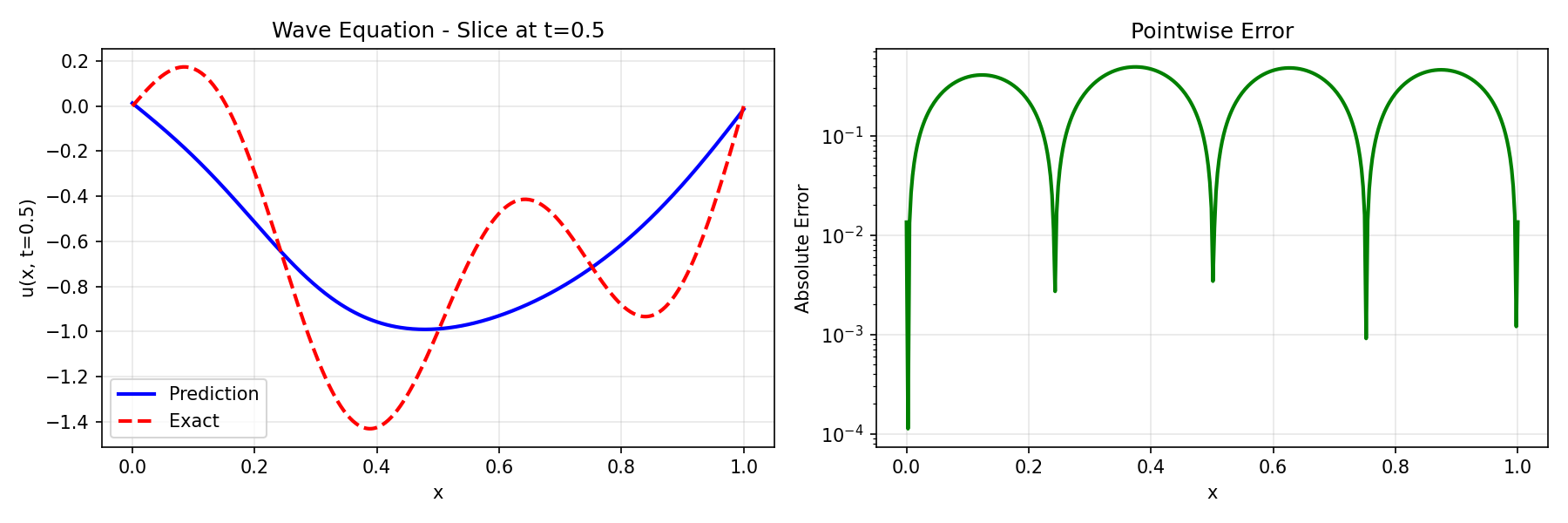}
    \caption{Baseline dual-network prediction at $t=0.5$ using tanh activations. Left: the predicted solution (blue) captures only the fundamental mode, missing the high-frequency oscillations present in the exact solution (red dashed). Right: pointwise absolute error on logarithmic scale.}
    \label{fig:wave_baseline}
\end{figure}

Figure~\ref{fig:wave_baseline} reveals the underlying cause of this limited performance: the predicted solution captures only the fundamental mode $\sin(\pi x)\cos(2\pi t)$, while the higher-frequency component $\frac{1}{2}\sin(4\pi x)\cos(8\pi t)$ is almost entirely absent. This behavior is consistent with the well-documented spectral bias phenomenon in neural networks with smooth activations~\cite{rahaman2019spectral, wang2021eigenvector, sitzmann2020implicit}, whereby standard architectures preferentially learn low-frequency functions and struggle to represent high-frequency components. 

\textbf{Two-nets and two-phase training (extended).} To address the spectral bias limitation, we introduce two extensions that preserve the core dual-network architecture. First, following Sitzmann et al.~\cite{sitzmann2020implicit}, we replace tanh with sinusoidal activations $\sigma(x) = \sin(\omega_0 x)$ (SIREN) in the domain network. Second, we add a modal prior regularization term that guides the learned Fourier coefficients toward their expected temporal evolution:
\begin{equation}
\mathcal{L}_{\text{modal}} = \sum_{n \in \mathcal{M}} \mathbb{E}_t\left[\left(a_n(t) - a_n(0)\cos(cn\pi t)\right)^2\right],
\end{equation}
where $a_n(t) = 2\int_0^1 u(x,t)\sin(n\pi x)\,dx$ and $\mathcal{M} = \{1, 4\}$ contains the dominant mode indices. The results in Table~\ref{tab:wave_extended} show that this extended configuration achieves mean MAE $(5.06 \pm 2.52) \times 10^{-2}$ and mean accuracy $88.02\%$, with the best run (seed~48) reaching MAE $3.0 \times 10^{-2}$ and accuracy $92.83\%$.

\begin{table}[H]
\centering
\caption{Two-nets, two-phase PINN with SIREN and modal prior for the wave equation (extended). Performance metrics (mean $\pm$ standard deviation) across different seed values.}
\label{tab:wave_extended}
\begin{tabular}{@{}lcccc@{}}
\toprule
& \multicolumn{3}{c}{Interior} & \\
Seed & MAE & RMSE & Rel.~$L^2$ & Accuracy (\%) \\
\midrule
42 & $3.5\times 10^{-2}$ & $4.5\times 10^{-2}$ & $8.1\times 10^{-2}$ & 91.86 \\
44 & $4.2\times 10^{-2}$ & $5.5\times 10^{-2}$ & $9.9\times 10^{-2}$ & 90.10 \\
46 & $9.3\times 10^{-2}$ & $1.2\times 10^{-1}$ & $2.2\times 10^{-1}$ & 77.53 \\
48 & $3.0\times 10^{-2}$ & $4.0\times 10^{-2}$ & $7.1\times 10^{-2}$ & 92.83 \\
50 & $5.3\times 10^{-2}$ & $6.8\times 10^{-2}$ & $1.2\times 10^{-1}$ & 87.80 \\
\midrule
\textbf{Mean $\pm$ Std} & $(5.06 \pm 2.52)\times 10^{-2}$ & $(6.56 \pm 3.22)\times 10^{-2}$ & $(1.18 \pm 0.60)\times 10^{-1}$ & $88.02 \pm 6.17$ \\
\bottomrule
\end{tabular}
\end{table}

\begin{figure}[H]
    \centering
     \begin{subfigure}[b]{0.45\textwidth}
         \centering
         \includegraphics[height=4.2cm,width=\textwidth]{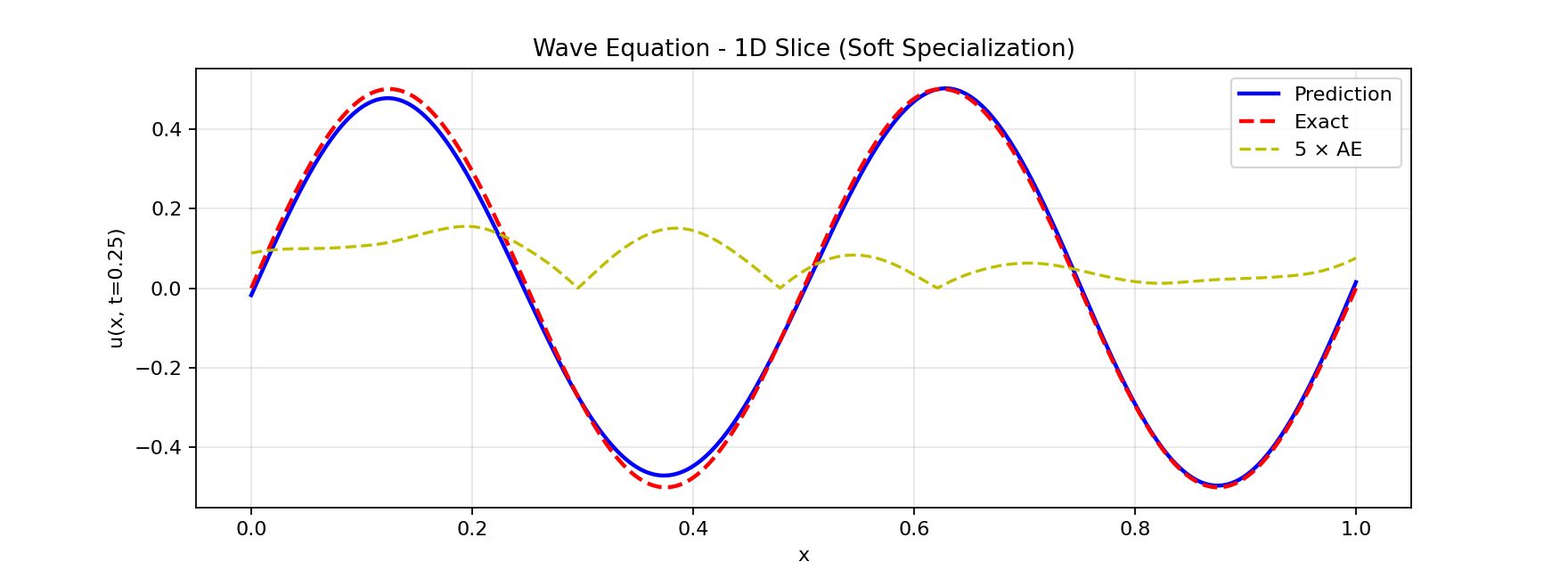}
         \caption{Comparison at $t=0.5$}
     \end{subfigure}
     ~~
    \begin{subfigure}[b]{0.45\textwidth}
    \centering
    \includegraphics[height=4.2cm,width=\textwidth]{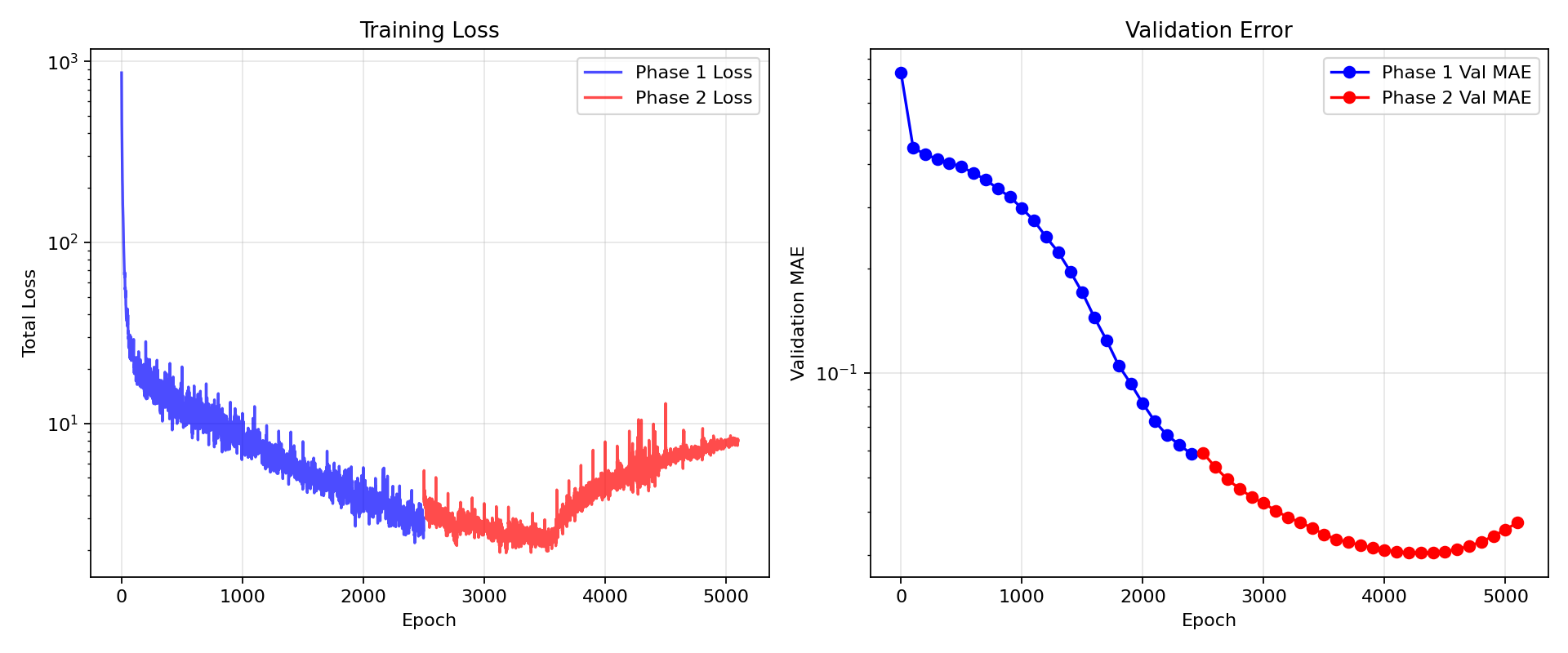}
    \caption{Training loss and validation error}
     \end{subfigure}
     ~~
     \begin{subfigure}[b]{\textwidth}
         \centering
         \includegraphics[height=4.2cm,width=\textwidth]{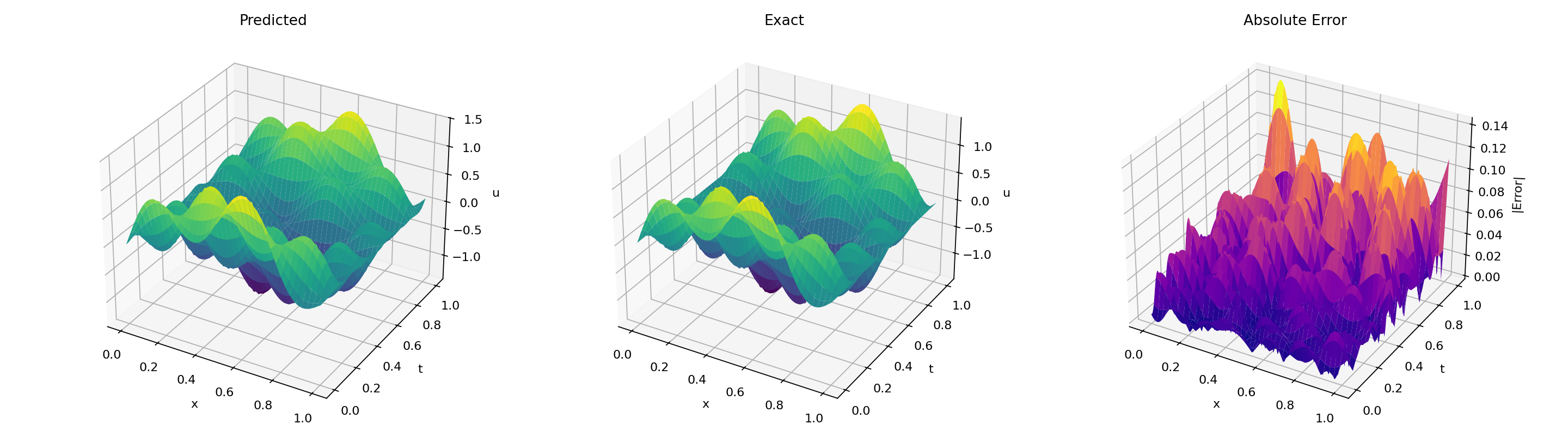}
         \caption{Surface plots}
     \end{subfigure}
    \caption{Best run results of wave equation with extended framework (seed~48).}
    \label{fig:wave_extended}
\end{figure}

Figure~\ref{fig:wave_extended} visualizes the best run (seed~48). Unlike the baseline, the extended model accurately reproduces both the fundamental and higher-frequency modes across the full $(x,t)$ domain.

The comparison in Table~\ref{tab:wave_comparison} confirms that the extended dual-network configuration substantially outperforms the baseline for this multi-scale benchmark. The $3.7\times$ reduction in mean MAE and $6\times$ reduction in best-case MAE demonstrate the effectiveness of SIREN activations and modal prior regularization for problems with high-frequency content. Importantly, all foundational elements of the framework remain intact: the shared physics residual, soft boundary-interior specialization, ALM-based enforcement, and the two-phase training curriculum.

\begin{table}[H]
\centering
\caption{Comparison between baseline and extended dual-network configurations for the wave equation.}
\label{tab:wave_comparison}
\begin{tabular}{lccc}
\hline
Method & MAE (interior) & Rel.~$L^2$ & Accuracy (\%) \\
\hline
Dual-net baseline (tanh) & $1.88 \times 10^{-1}$ & $4.15 \times 10^{-1}$ & 58.5 \\
Dual-net extended, mean & $5.06 \times 10^{-2}$ & $1.18 \times 10^{-1}$ & 88.0 \\
Dual-net extended, best (seed 48) & $3.0 \times 10^{-2}$ & $7.1 \times 10^{-2}$ & 92.8 \\
\hline
\end{tabular}
\end{table}

\section{Conclusion}
We introduced a dual-network PINN architecture that decomposes the solution into interior and boundary components, coupling them through a unified physics loss, distance-weighted specialization, and an augmented Lagrangian formulation for boundary conditions. This design addresses the long-standing challenges of single-network PINNs on multi-scale PDEs, namely optimization interference between interior physics and boundary enforcement, difficulty in resolving sharp gradients, and sensitivity to penalty hyperparameters. By using soft role priors that are cosine-annealed during training and a two-phase sampling curriculum that first stabilizes boundary satisfaction and then concentrates resolution where it is most needed, our method reliably captures both global solution structure and localized features.
Across four benchmark problems, including Laplace equation, Poisson equation, and 1D Fokker-Planck equations, the proposed framework consistently outperforms standard PINNs, yielding substantially lower relative errors, tighter boundary adherence, and improved resolution of steep solution variations. Ablation studies further confirm the complementary roles of boundary-interior specialization, annealed role regularization, and curriculum sampling. In general, the proposed approach is simple to implement, incurs minimal computational overhead, and integrates seamlessly with existing PINNs workflows. We believe that this framework offers a practical step toward making physics-informed models more robust for PDEs with sharp transitions and complex boundary conditions, and opens avenues for extensions to higher-dimensional, time-dependent, and stochastic systems.
\\
\\
\textbf{Acknowledgements}
\\
{ This publication was partially funded by the PhD program in Mathematics and Computer Science at University of Calabria, Cycle XXXIX with the support of a scholarship financed by DM
118/2023 (CUP H23C23001760005), based on the PNRR funded by the European Union.
}
\\
\\
\textbf{Funding} Open access funding provided by Università della Calabria within the CRUI-CARE Agreement.
\\
\\
\textbf{Data Availability}
\\
{{All the data generated or analyzed during this study are included in this published article.}}
\\
\section*{Declarations}

\textbf{Conflicts of Interest}
{The authors declare that they have no conflict of interest.}
\\

\bibliography{library}
\bibliographystyle{unsrt} 

\end{document}